\documentclass[preprint]{elsarticle}
\usepackage{amssymb}
\usepackage{amsmath}
\usepackage{latexsym}
\usepackage{subcaption}
\usepackage{url}

\usepackage{diagbox}

\usepackage{algorithm}
\usepackage{algpseudocode}
\usepackage{algorithmicx}
\usepackage{tikz}
\usetikzlibrary{tikzmark,decorations.pathreplacing}

\journal{Journal of Computational Physics}

\begin{document}

\begin{frontmatter}
\title{A parallel-in-time approach for accelerating direct-adjoint studies}%
\author[1]{C. S. Skene\corref{cor1}\fnref{UCLA}}
\ead{csskene@ucla.edu}
\author[1]{M. F. Eggl \fnref{Princeton}}
\author[1]{P. J. Schmid}
\address[1]{Department of Mathematics, Imperial College London, London SW7 2AZ, United Kingdom}
\fntext[UCLA]{Present address: Department of Mechanical and Aerospace Engineering, University of California, Los Angeles, CA 90095, USA}
\fntext[Princeton]{Present address: Department of Mechanical and Aerospace Engineering, Princeton University, NJ 08544, USA}
\cortext[cor1]{Corresponding author}

\begin{abstract}
  Parallel-in-time methods are developed to accelerate the
  direct-adjoint looping procedure. Particularly, we utilize the
  {\tt{Paraexp}} algorithm, previously developed to integrate
  equations forward in time, to accelerate the direct-adjoint looping
  that arises from gradient-based optimization. We consider both
  linear and non-linear governing equations and exploit the linear,
  time-varying nature of the adjoint equations. Gains in efficiency
  are seen across all cases, showing that a Paraexp based parallel-in-time approach is feasible for the acceleration of direct-adjoint studies. This
  signifies a possible approach to further increase the run-time
  performance for optimization studies that either cannot be
  parallelized in space or are at their limit of efficiency gains for
  a parallel-in-space approach.
\end{abstract}


\end{frontmatter}

\section{Introduction}
The use of direct-adjoint looping techniques to compute the extremum
of an objective functional, subject to a state variable satisfying
complementary governing (direct) equations, is the subject of a wide
range of optimization studies. Direct-adjoint looping is a technique
that originated in control theory and was pioneered in fluid dynamics
by \citet{Jameson1988} to find the optimal aerofoil shape that
minimizes drag. These Lagrangian optimization problems are formulated
with the governing equations included as constraints. The direct
(governing) equations are integrated forward in time yielding a
solution for a set of control parameters. Afterwards, the adjoint
equations are integrated backwards in time, providing the gradient of
the objective functional with respect to the parameters. This
information can then be used as part of a gradient-based optimization
routine to find the extremum of the objective functional via repeated
applications of this direct-adjoint loop. Recent examples include,
among others, the optimization of mixing in binary
fluids~\citep{Foures2014,Marcotte2018,Eggl2018,eggl2019,Eggl2019b},
finding the minimal seed that triggers turbulence in pipe
flow~\citep{KerswellTurb} and determining the optimal place to ignite
a diffusion flame~\citep{ubaidFlame}.

While the inclusion of the adjoint procedure significantly speeds up
the calculation of the gradient information, several complications can
arise due to the use of the adjoint. Firstly, the inherently iterative
nature of direct-adjoint looping can cause optimization studies
performed in this manner to become prohibitively slow, despite the
speedup achieved through the efficient computation of the
gradient. Another complication arises when the adjoint equations
depend on the direct solution. This necessarily implies that the full
direct solution must be saved and injected into the adjoint equation
at specific points in time to obtain accurate solutions for the
adjoint variables, leading to high memory costs. To circumvent these
memory costs the checkpointing library {\tt{revolve}}~\citep{revolve}
can be used. In this way, the direct solution is saved at various
checkpoints in time, and, at the expense of recomputing the solution,
the direct-adjoint loop can then be solved between checkpoints leading
to a reduced memory cost. These two complications, among many others,
can result in run-times that may be unfeasible for some optimization
studies. Traditionally, to alleviate these problems and retain viable
run-times, parallel-in-space approaches are used to speed up the
calculation, examples of which include work by
\citet{pekurovsky2012,Laizet2010,kallala2019}. However, when
encountering problems that cannot be spatially parallelized or where
the parallelization in space has reached its maximum efficiency, gains
in run-time efficiency must be achieved through other means. In our
study, we consider one such alternative method to accelerate scenarios
where more computational power is available: parallelization in time.

A wide variety of parallel-in-time algorithms have been developed for parallelizing evolution equations, and for an in-depth review we refer the reader to the review papers of \citet{gander50years} and \citet{Ong2020}. These methods aim to solve an equation more efficiently in parallel through a decomposition of the time domain. As the solution of an evolution equation at a future time depends on its full previous time history, parallelizing-in-time such processes is non-trivial and several approaches have arisen to tackle this challenge. The first successful attempt to overcome this difficulty was presented in the work of \citet{Nievergelt1964} who decomposed the full time domain into subdomains and calculated the solution in parallel using a shooting method. This algorithm is direct, however, since then other shooting methods that are iterative have been developed, notably, the widely studied {\tt Parareal} algorithm \citep{lions2001} on which several other algorithms, such as parallel-implicit-time-integrator ({\tt PITA}) \citep{Farhat2003}, are based. As well as shooting based methods other techniques which decompose time differently exist; decomposing space-time into subdomains in space that span the whole time domain gives rise to waveform relaxation methods \citep{Lelarasmee1982,Gander1996,Giladi2002,Kwok2014,Mandal2014}, whereas parallelizing the solver over the whole of space-time leads to multigrid methods \citep{HortonVandewalle1995,EmmettMinion2012,Falgout2014,Gander2016}. Further approaches include algorithms based on diagonalizing the time-stepping matrix so that all the time-steps can be solved in parallel~\citep{MADAY2008113} and are known as {\tt ParaDiag} algorithms (see~\citep{GanderEtAl2020} for further information).

In our study we consider the {\tt{Paraexp}} algorithm developed
by~\citet{GanderParaExp} which is a direct method particularly suited for hyperbolic and diffusive problems \citep{gander50years}. This algorithm achieves parallelization-in-time for a linear equation by isolating the
homogeneous and inhomogeneous parts. The time domain is then
partitioned and the inhomogeneous equations are solved on each time
partition with a zero initial condition. Solving the homogeneous
equations allows for propagation of the correct initial conditions
throughout the domain, with the full solution obtained as a
superposition of all inhomogeneous and homogeneous parts. Although the
homogeneous equations can span more than one time partition, and hence more initial value problems are solved than the original problem consisted of, the fact that the homogeneous equations can often be efficiently solved via exponential time integration, examples of which can be seen in~\citet{schulze2009} and~\citet{realLeja}, leads to an overall speedup. Furthermore, this algorithm can be extended to non-linear equations via an iterative approach~\citep{ParaExpNL,Gander2017ANP}.

Parallel-in-time methods have proved to be a useful choice for parallelising direct studies, especially when spatial parallelization has saturated. For example, the work of \citet{Farhat2003} develops and applies {\tt PITA} to fluid, structure and fluid-structure interaction problems, and the work of \citet{Guillaume2002} considers a non-linear extension of the {\tt Parareal} algorithm and implements it on the pricing of an American put option using the Black-Scholes equation. Further to these, recently the applicability of the {\tt Parareal} algorithm for accelerating the time integration of kinetic dynamos \citep{Clarke2020a} and Rayleigh-B{\'{e}}nard convection \citep{Clarke2020b} has been considered. Aside from being used to accelerate direct studies, parallel-in-time methods have also been extended to optimization studies. In the work of \citet{Gunther2018} and \citet{Gunther2019} the {\tt XBraid} library, which utilizes a multigrid reduction-in-time technique \citep{FriedhoffEtAl2013}, is extended to accelerate optimization studies. Likewise, the {\tt PFASST} algorithm \cite{EmmettMinion2012} has also been used for PDE optimization \citep{Gotschel2018,Gotschel2019}. In addition to these multigrid-based approaches, several optimization algorithms based on {\tt Parareal} have been introduced such as in the work of~\citet{Maday2002},~\citet{Maday2003} and more recently in the algorithm {\tt ParaOpt}~\citep{Gander2020}.

This present work aims to accelerate direct-adjoint looping by
extending the {\tt{Paraexp}} algorithm to also include the adjoint
component. We consider both linear and non-linear equations, as well
as the possible inclusion of a checkpointing scheme. For the case of
non-linear governing equations two algorithms are developed; one in
which the direct equation is solved iteratively in parallel, and
another in which the direct equation is solved in series. In both
cases, the adjoint equation is solved in parallel. When the direct
equation is solved in series we obtain a speedup by overlapping the
direct and adjoint solutions using a particularly efficient time
partition. For all algorithms we derive theoretical expressions for
the speedups. Furthermore, we describe in a step-by-step fashion the
developed parallel-in-time algorithms and apply them to the 2D
advection-diffusion equation and viscous Burgers' equation for the
linear and nonlinear cases, respectively. In all cases we observe a
marked improvement in run-time efficiency in line with the theoretical
predictions.

The rest of this article is organized in the following manner; \S2
gives a brief overview of adjoint-looping techniques before moving on
to \S3 where parallel-in-time algorithms are presented for the direct
and adjoint system. Next, \S4 covers the numerical implementation, and
the resulting numerical experiments are compared to theoretical
scalings. Finally, conclusions are offered in section \S5.

\section{Direct-adjoint looping}\label{sec:directadjoint}

We begin by briefly outlining a general direct-adjoint problem
\citep{Jameson1988}. Consider the non-linear governing (direct)
equation

\begin{equation}
  \frac{\partial {\bf q}}{\partial t}=\boldsymbol{\mathcal{N}}({\bf
    q},{\bf p}), \qquad {\bf q}(0)={\bf q}_0,
    \label{equ:state}
\end{equation}
where $\boldsymbol{\mathcal{N}}({\bf q},{\bf p})$ represents the
right-hand side of an equation dependent upon the state ${\bf q}$ and
parameters ${\bf p}$. We seek to minimize an objective functional
$\mathcal{J}$ that also can depend on the state ${\bf q}$ and
parameters ${\bf p}$. To this end, we form the Lagrangian
\begin{equation}
    \mathcal{L}({\bf q},{\bf p})=\mathcal{J}({\bf q},{\bf p})-\left(\left \langle {\bf q}^\dagger,
    \frac{\partial {\bf q}}{\partial t}-{ \bf\mathcal{N}}({\bf q},{\bf
      p})\right\rangle + \textrm{c.c.}\right),
\end{equation}
where the Lagrange multiplier ${\bf q}^\dagger$ is the adjoint
variable. Note that we have used the abbreviation c.c. to denote the complex conjugate, which is included to ensure the Lagrangian is real valued. We take the inner products and the cost functional to be of
the form
\begin{eqnarray}
  \left \langle {\bf a},{\bf b}\right \rangle &=&\int_0^T [ {\bf
      a},{\bf b}]\; \textrm{d}t,\\
  \left [ {\bf a},{\bf b}\right ] &=&
  \begin{cases}
    {\bf a}^H{\bf b}, \;& \textrm{discrete case,} \\
    \int_\Omega{\bf a}^H{\bf b}\;\textrm{d}V &\textrm{continuous case,}
  \end{cases}\\
  \mathcal{J}(\bf{q},{\bf p},{\bf q}(T))&=&
  \begin{cases}
    \int_0^T \mathcal{J}_I({\bf q},{\bf p}) \,
    \textrm{d}t+\mathcal{J}_T({\bf p},{\bf q}(T)) &\textrm{discrete
      case,} \\
    \int_0^T\int_\Omega\mathcal{J}_I({\bf q},{\bf p}) \,
    \textrm{d}V\textrm{d}t+\mathcal{J}_T({\bf p},{\bf q}(T))
    &\textrm{continuous case,}
  \end{cases}
\end{eqnarray}
where differing inner products and cost functionals are introduced for
the discrete and continuous cases, respectively. In our formulation,
$\Omega$ represents the spatial domain, and the cost functional is
split into an integral part, $\mathcal{J}_I$, and a part that solely
depends on the final time, $\mathcal{J}_T$.

To minimize our objective functional subject to the constraint that
our state satisfies the governing equation~(\ref{equ:state}), we must
ensure that all first variations of the Lagrangian are zero. Taking
the variation with respect to the conjugate adjoint variable enforces
the state equation~(\ref{equ:state}). Before we take the variations
with respect to the adjoint variable ${\bf q}^\dagger$, we utilize
integration by parts in time to rewrite the Lagrangian as follows

\begin{equation}
  \mathcal{L}=\mathcal{J}+\left(-[{\bf q}^\dagger,{\bf q}]_0^T+\left
  \langle \frac{\partial{\bf q}^\dagger}{\partial t},{\bf q} \right
  \rangle + \left \langle {\bf q}^\dagger,{ \bf\mathcal{N}}({\bf
    q},{\bf p})\right\rangle+ \textrm{c.c.}\right).
\end{equation}
Upon taking the variations of the Lagrangian in this form with respect
to the state ${\bf q}$, we arrive at the adjoint equation
\begin{equation}
  \frac{\partial {\bf q}^\dagger}{\partial t}= -\frac{\partial{
      \bf\mathcal{N}}({\bf q},{\bf p})}{\partial {\bf q}}^\dagger{\bf
    q}^\dagger-\overline{{\frac{\partial \mathcal{J}_I}{\partial {\bf
          q}}}}
    \label{equ:adjEquation}
\end{equation}
where $\overline{(\cdot)}$ denotes the complex conjugate. The adjoint
operator must be found such that the adjoint relation
\begin{equation}
  \left[{\bf a},\frac{\partial{ \bf\mathcal{N}}({\bf q},{\bf p})}{\partial {\bf q}}{\bf b}\right] = \left[\frac{\partial{
        \bf\mathcal{N}}({\bf q},{\bf p})}{\partial {\bf
        q}}^\dagger{\bf a},{\bf b}\right]
\end{equation}
is satisfied for all vectors $\bf{a}$ and ${\bf b}$. In the continuous
case, this relation can be used to find the adjoint operator via
integration by parts in space. For the discrete case, on the other
hand, simple linear algebra manipulations can be used to find the
adjoint operator
\begin{equation}
  \frac{\partial{ \bf\mathcal{N}}({\bf q},{\bf p})}{\partial {\bf
      q}}^\dagger=\frac{\partial{ \bf\mathcal{N}}({\bf q},{\bf p})}{\partial {\bf q}}^H.
\end{equation}

As the adjoint equation~(\ref{equ:adjEquation}) must be integrated
backwards in time from the final time $T$, we require a final time
condition for the adjoint variable, ${\bf q}^\dagger$. This can found
by ensuring the first variation with respect to ${\bf q}(T)$ is zero,
i.e.,
\begin{equation}
  {\bf q}^\dagger(T)=\overline{\frac{\partial \mathcal{J}_T}{\partial{\bf q}(T)}}
  \label{equ:adjFinalTime}.
\end{equation}
Finally, it remains to take the first variation with respect to ${\bf p}$. In general, unless ${\bf p}$ is a local minimum of the
Lagrangian, the first variation with respect to ${\bf p}$ will be
non-zero. Instead, we obtain the gradient of the Lagrangian as
\begin{equation}
  \frac{\partial \mathcal{L}}{\partial {\bf p}}= \frac{\partial
    \mathcal{J}}{\partial {\bf p}}+\left\langle {\bf
    q}^\dagger,\frac{\partial{ \bf\mathcal{N}}({\bf q},{\bf p})}{\partial {\bf p}}\right \rangle.
    \label{equ:paramupdate}
\end{equation}
Likewise, if the cost functional is to be optimized with respect to
the initial condition, the associated gradient becomes
\begin{equation}
  \frac{\partial \mathcal{L}}{\partial{\bf q}_0}= \frac{\partial
    \mathcal{J}}{\partial{\bf q}_0}+\overline{{\bf q}^\dagger(0)}.
    \label{equ:initupdate}
\end{equation}

In practice, achieving an analytic expressions for all required
quantities is not feasible. Therefore, the governing equations and
adjoint equations are solved explicitly, while some initial guess is
employed for our control variable, ${\bf p}$ in the case of forcing or
${\bf q_0}$ in the case of optimal initial condition. An iterative
procedure, in conjunction with either equation~(\ref{equ:paramupdate})
or~(\ref{equ:initupdate}) and a gradient-descent routine, is then
utilized to drive the control to a value that achieves a local minimum
for the chosen cost functional.

\section{Parallel-in-time algorithms}\label{sec:algorithms}

\subsection{Linear governing equation}\label{sec:linearAlg}

We now turn to the discussion of the parallel-in-time approach,
particularly the {\tt{Paraexp}} algorithm developed by
\citet{GanderParaExp}. The fundamental aspect of this approach is that
a linear equation can be split into an inhomogeneous and homogeneous
component. Consider the general linear equation
\begin{equation}
  \frac{\partial {\bf q}}{\partial t}={\bf A}{\bf q}+{\bf f},
  \qquad{\bf q}(0)={\bf q}_0,
  \label{equ:lin}
\end{equation}
where ${\bf q}\in\mathbb{C}^n$, ${\bf A}\in\mathbb{C}^{n\times n}$ and the inhomogeneity arises from the source term
${\bf f}\in\mathbb{C}^n$. Following \citet{ParaExpNL}, we make the substitution
$\tilde{\bf{q}}={\bf q}-{\bf q}_0$ to obtain the equation
\begin{equation}
  \frac{\partial \tilde{\bf q}}{\partial t}={\bf A}\tilde{\bf q}+{\bf A}\tilde{\bf
    q}_0+{\bf f}, \qquad \tilde{\bf{q}}(0)={\bf 0},
  \label{equ:linzero}
\end{equation}
which absorbs the initial condition into the governing equation.

To solve this equation over the range $t\in [0,T]$ in parallel, we
employ the {\tt{Paraexp}} algorithm. For $N$ processors the time
domain is decomposed into the partition $0=T_0<T_1<...<T_{N}=T$. The
$p^\textrm{th}$ processor, where available processors are labelled
$p\in \{0,1,...,N-1\}$, then solves the inhomogeneous equation
\begin{equation}
  \frac{\partial {\bf q}_{I,p}}{\partial t}={\bf A}{\bf q}_{I,p}+{\bf
    A}{\bf q}_0+{\bf f}, \qquad {\bf q}_{I,p}(T_p)={\bf 0},
\label{equ:linInhomEqus}
\end{equation}
on $[T_{p},T_{p+1}]$. Thus the inhomogeneous part of the equation is
solved on each time interval with a zero initial condition. Following
this we solve the homogeneous problems
\begin{equation}
  \frac{\partial {\bf q}_{H,j}}{\partial t}={\bf A}{\bf q}_{H,j},
  \qquad {\bf q}_{H,j}(T_j)={\bf q}_{I,j-1}(T_j),
  \label{equ:hom}
\end{equation}
on the intervals $[T_{j},T]$ for $j\in\{1,...N-1\}$.  At first it may
seem that solving the equation in series would be faster than the
parallel implementation, as the homogeneous components require solving
over at least the interval $[T_1,T_N]$, which in addition to the
inhomogeneous evaluation covers the full domain. However, the speedup
is achieved by recognizing that the homogeneous equations have the
ability to be solved far more efficiently than the inhomogeneous
equation, for example by using exponential
time-stepping~\citep{krylovApprox}. Hence, the slow-to-solve
inhomogeneous parts are solved in parallel first and then the
homogeneous parts are solved rapidly to propagate the correct initial
conditions throughout the time domain. The full solution on the
interval $[T_p,T_{p+1}]$ can then be found by the sum
\begin{equation}
  {\bf q} = {\bf q}_0+{\bf q}_{I,p} + \sum_{j=1}^{p} {\bf q}_{H,j}.
  \label{equ:fullsol0}
\end{equation}

\begin{figure}
  \centering
  \begin{tabular}{cc}
    \begin{subfigure}[t]{0.4\textwidth}
      \includegraphics[width=\linewidth]{ 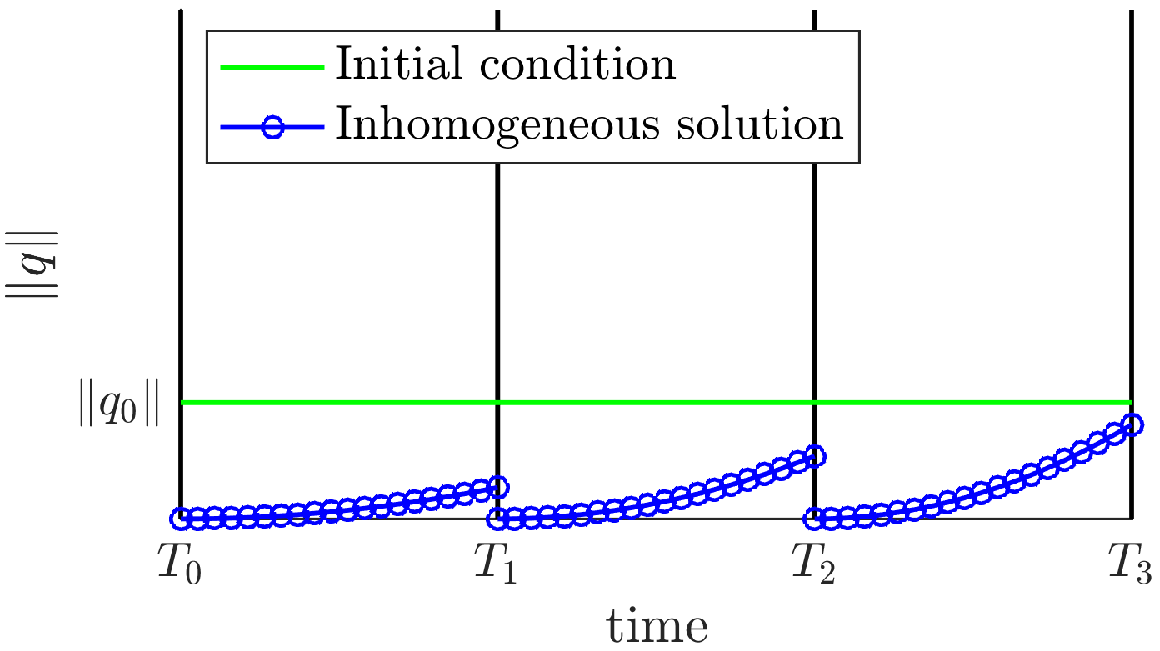}
      \caption{First, the inhomogeneous equations are solved.}
      \label{fig:Inhom1}
    \end{subfigure}
    & \hspace{1truecm}
    \begin{subfigure}[t]{0.4\textwidth}
      \includegraphics[width=\linewidth]{ 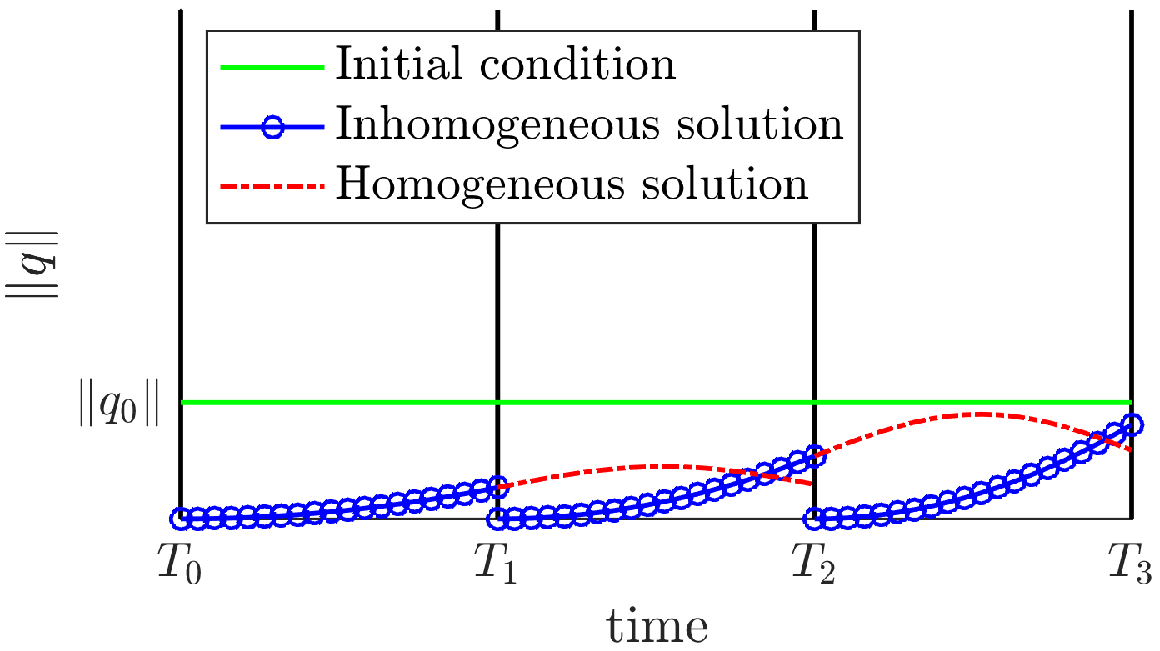}
      \caption{Next, the homogeneous parts are solved on each processor.}
      \label{fig:hom2}
    \end{subfigure} \\
    \begin{subfigure}[t]{0.4\textwidth}
      \includegraphics[width=\linewidth]{ 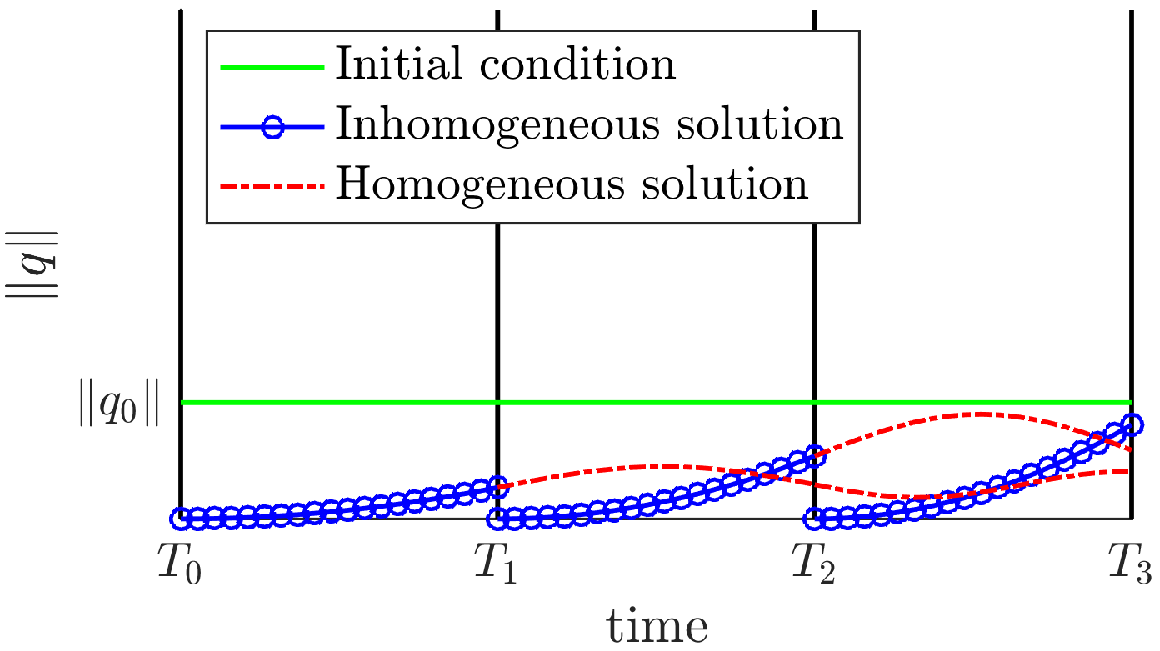}
      \caption{The homogeneous solution on the second processor is
        then solved to the final time on the third processor.}
      \label{fig:hom3}
    \end{subfigure}
    & \hspace{1truecm}
    \begin{subfigure}[t]{0.4\textwidth}
      \includegraphics[width=\linewidth]{ 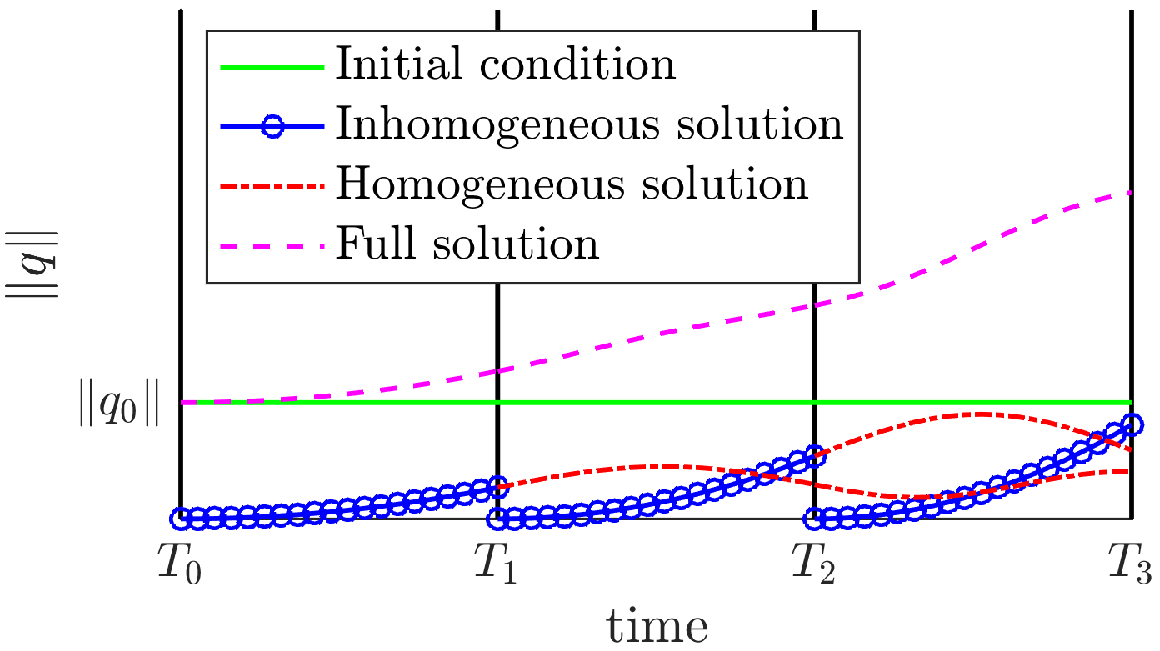}
      \caption{Finally, we sum all parts on each processor to find the full solution.}
      \label{fig:sol4}
    \end{subfigure}
  \end{tabular}
  \caption{An example of an equation being solved in parallel with
    three processors. The subfigures show the progression of the
    algorithm.}
  \label{fig:parallelGuide}
\end{figure}

In their implementation of the {\tt{Paraexp}} algorithm,
\citet{GanderParaExp} suggest solving problems ${\bf q}_{I,p}$ and
${\bf q}_{H,p+1}$ on the same processor so that communication between
processors is only necessary when computing the
sum~(\ref{equ:fullsol}) at the end. Here, however, we will consider keeping the
time domain completely separate, i.e., the $p^\textrm{th}$ processor
will always solve equations on its time partition
$[T_p,T_{p+1}]$. Therefore, whilst the inhomogeneous problems are
separated by processors, an individual solution for the homogeneous
part ${\bf q}_{H,j}$ will span multiple processors, resulting in
necessary communications between adjacent processors. The reason for
this particular arrangement is due to the inclusion of the adjoint
solver. The adjoint will require full state information on each time
interval and therefore it is vital for efficiency that this is
available without communicating a large amount of data across
processors. Hence, we rewrite the homogeneous problems~(\ref{equ:hom})
as

\begin{equation}
  \frac{\partial {\bf q}_{H,p,j}}{\partial t}={\bf A}{\bf q}_{H,p,j},
  \label{equ:linHomEqus}
\end{equation}
for $p>0$ and $j\in\{0,...,p-1\}$, to be solved on $[T_p,T_{p+1}]$
with the initial conditions
\begin{equation}
  {\bf q}_{H,p,j}(T_p)  =
  \begin{cases}
    {\bf q}_{I,p-1}(T_p),  \; & j=0, p \neq 0,\\
    {\bf q}_{H,p-1,j-1}(T_p),  \; & 1\leq j\leq p-1, p \neq 0,1.
  \end{cases}
\end{equation}
Here the subscript $p$ completely defines the processor that the
relevant equation is solved on, and the $p^\textrm{th}$ processor
solves $p$ homogeneous problems. We note that the subscript $j$ is not
the same as in equation~(\ref{equ:hom}) but is instead used to denote
the different homogeneous problems that must be solved for each
processor. Each processor (except the zero$^{\text{th}}$) must solve
the $j=0$ homogeneous problem that stems from the final time
inhomogeneous state from the previous processor. Then processor $p$
must solve the additional homogeneous problems $j=1,...,(p-1)$ that
arise from the homogeneous problems started on previous processors
that must be integrated to the final time. By solving the equations
this way the full solution on the interval $[T_p,T_{p+1}]$ becomes
\begin{equation}
  {\bf q} = {\bf q}_0+{\bf q}_{I,p} + \sum_{j=0,p\neq 0}^{p-1} {\bf q}_{H,p,j}.
  \label{equ:fullsol}
\end{equation}

It should be noted that solving the homogeneous equations in this way may introduce some lag that is not present in the {\tt{Paraexp}} algorithm considered by \citet{GanderParaExp}. This lag occurs when the homogeneous equations do not take the same time to be solved on each processor, meaning that some processors must wait to send their final states to the next processor to continue the integration. However, this introduced lag is assumed to be smaller than the cost of communicating the whole direct solution across all processors in order to solve the adjoint equations. In cases where the efficiency deterioration due to this lag becomes comparable to, or larger than, the cost of distributing the direct solution, the homogeneous solutions should be computed in the manner suggested by \citet{GanderParaExp}, and this is the approach considered in \ref{sec:hyb}.

Figure~\ref{fig:parallelGuide} illustrates the progression of the
algorithm for a three processor arrangement. Firstly, the
inhomogeneous equations~(\ref{equ:linInhomEqus}) are solved as
displayed in figure~\ref{fig:Inhom1}. Figures~\ref{fig:hom2}
and~\ref{fig:hom3} then show the homogeneous
equations~(\ref{equ:linHomEqus}) being solved. We see the $j=0$
homogeneous states being solved on all processors in
figure~\ref{fig:hom2}, with the final time solution from processor
$p=1$ being passed to processor $p=2$ to solve the homogeneous state
$j=1$ in figure~\ref{fig:hom2}. Finally, the sum~(\ref{equ:fullsol0})
is computed on each processor to obtain the full solution.

Now that the algorithm for the direct equation has been presented we
turn to the discussion on the inclusion of the adjoint solver. The
adjoint equation to~(\ref{equ:lin}) for a given cost functional (see
section~\ref{sec:directadjoint}) is
\begin{equation}
  \frac{\partial \tilde{\bf q}^{\dagger} }{\partial t}=-{\bf A}^H\tilde{\bf
    q}^{\dagger} -{\bf A}^H{\bf
    q}^\dagger(T)-\overline{\frac{\partial \mathcal{J}_I}{\partial
      {\bf q}}},\qquad \;\tilde{\bf q}^{\dagger} (T)={\bf 0},
\end{equation}
where we now have to integrate backwards in time, and the final time
condition was absorbed into the equation using the substitution $\tilde{\bf
  q}^{\dagger}={\bf q}^\dagger-{\bf q}^\dagger(T)$. The linear nature
of this equation means that we can use the {\tt{Paraexp}} algorithm in
exactly the same manner as previously described. The adjoint
inhomogeneous problems that must be solved are
\begin{equation}
  \frac{\partial {\bf q}^\dagger_{I,p}}{\partial t}=-{\bf A}^H{\bf
    q}^\dagger_{I,p}-{\bf A}^H{\bf
    q}^\dagger(T)-\overline{\frac{\partial \mathcal{J}_I}{\partial
      {\bf q}}}, \qquad {\bf q}^\dagger_{I,p}(T_{p+1})={\bf 0},
  \label{equ:adjInhom}
\end{equation}
on $[T_{p+1},T_{p}]$. Similarly to the direct case, the adjoint
homogeneous problems are
\begin{equation}
  \frac{\partial {\bf q}^\dagger_{H,j}}{\partial t}=-{\bf A}^H{\bf
    q}^\dagger_{H,j}, \qquad {\bf q}^\dagger_{H,j}(T_{j+1})={\bf
    q}^\dagger_{I,j+1}(T_{j+1}),
\label{equ:adjHomNormal}
\end{equation}
on $[T_{j+1},0]$ and for $j\in\{N-2,...0\}$. Once again, we keep the
adjoint homogeneous problems local to each processor's time partition
and instead solve the problems
\begin{equation}
  \frac{\partial {\bf q}^\dagger_{H,p,j}}{\partial t}=-{\bf A}^H{\bf
    q}^\dagger_{H,p,j},
  \label{equ:adjHomLocal}
\end{equation}
on $[T_{p+1},T_{p}]$ with the initial conditions
\begin{equation}
  {\bf q}^\dagger_{H,p,j}(T_{p+1})  =
  \begin{cases}
    {\bf q}^\dagger_{I,p+1}(T_{p+1}), \; & j=0, p\neq N-1,\\
    {\bf q}^\dagger_{H,p+1,j-1}(T_{p+1}), \; & 1\leq j \leq (N-1)-p, \;\;
    p\neq N-1,N-2.
  \end{cases}
\end{equation}
This time the $p^\textrm{th}$ processor solves $(N-1)-p$ homogeneous
equations. Note that these equations can be solved forwards in time by
making the substitution $\tau=T_{p}+T_{p+1}-t$ on each time
interval. The full adjoint solution on processor $p$ is then given as
the sum
\begin{equation}
  {\bf q}^\dagger = {\bf q}^\dagger(T)+{\bf q}^\dagger_{I,p} +
  \sum_{j=0,p\neq N-1}^{(N-1)-p-1} {\bf q}^\dagger_{H,p,j},
  \label{equ:fulladjsol}
\end{equation}
on $[T_p,T_{p+1}]$.

\subsection{Non-linear governing equation solved in parallel}
\label{sec:nonLinParallel}

We now extend the parallel-in-time techniques to a direct-adjoint loop
for which the governing equation is non-linear. As the adjoint
equation is always a linear (and possibly linear, time-varying)
equation, the {\tt{Paraexp}} algorithm as presented in
section~\ref{sec:linearAlg} can always be used to solve the adjoint
equation in parallel. To obtain speedup for the direct equation we
present two approaches. The first uses a non-linear extension of the
{\tt{Paraexp}} algorithm, namely the method of~\citet{ParaExpNL}, to
accelerate the direct equation. Note that other non-linear extensions
of the {\tt{Paraexp}} algorithm are available, see for example the
proceedings of \citet{Gander2017ANP}, but we do not consider them
here. Specifically, the approaches of~\citet{Gander2017ANP}
and~\citet{ParaExpNL} are both iterative approaches and therefore
should yield fairly similar scalings.

For the non-linear governing equation we take the general
equation~(\ref{equ:state}) from section~\ref{sec:directadjoint}. To
parallelize the time-integration of this equation we proceed using the
method of~\citet{ParaExpNL} by rewriting it in the form
\begin{equation}
  \frac{\partial {\bf q}}{\partial t}={\bf A}{\bf q}+{\bf N}({\bf
    q})+{\bf f}, \; {\bf q}(0)={\bf q}_0,
  \label{equ:nl} 
\end{equation}
where we have split the operator $\boldsymbol{\mathcal{N}}$ into its
linear part ${\bf A}$ and genuine non-linear part ${\bf N}({\bf
  q})$. We then introduce the following recurrence relation
\begin{equation}
  \frac{\partial {\bf q}_{i+1}}{\partial t}={\bf A}{\bf q}_{i+1}+{\bf
    N}({\bf q}_i)+{\bf f}, \qquad {\bf q}_{i+1}(0)={\bf q}_0,
  \label{equ:nlIter}
\end{equation}
and note that if we achieve convergence for our series of ${\bf q}_i$,
i.e. ${\bf q}_{i+1}\approx{\bf q}_{i}$, then~(\ref{equ:nlIter}) is
equivalent to equation~(\ref{equ:nl}). The benefit of this altered
equation is that this equation is linear in ${\bf q}_{i+1}$, as the
non-linearity only acts on the previous iterate and becomes a source
term only.

To achieve the convergence necessary for our recurrence relation to
equate to the full equation, \citet{ParaExpNL} augment this equation
with
\begin{equation}
  \frac{\partial \tilde{\bf q}_{i+1}}{\partial t}=({\bf A}+{\bf J}_i)\tilde{\bf q}_{i+1}+{\bf A}{\bf q}_0+{\bf N}(\tilde{\bf q}_i)+{\bf f}-{\bf
    J}_i\tilde{\bf q}_i, \qquad \tilde{\bf q}_{i}(0)={\bf 0},
\end{equation}
where the average Jacobian 
\begin{equation}
  {\bf J}_i = \frac{1}{T}\int_0^T \left .\frac{\partial {\bf N}({\bf
      q})}{\partial {\bf q}}\right|_{{\bf q}={\bf q}_i} \;\textrm{d}t.
\end{equation}
has been added. Once again the substitution $\tilde{\bf q}={\bf q}-{\bf
  q}_0$ is used to absorb the initial condition into the equation.

The non-linear equations have now been reduced to a set of equivalent
linear ones and can be solved using the {\tt{Paraexp}} method
described in the previous section, giving the inhomogeneous problems
as
\begin{equation}
  \frac{\partial {\bf q}_{I,p,i+1}}{\partial t}=({\bf A}+{\bf
    J}_{i,p}){\bf q}_{I,p,i+1}+{\bf A}{\bf q}_0+{\bf N}({\bf
    q}_{I,p,i})+{\bf f}-{\bf J}_{i,p}{\bf q}_{I,p,i}, \; {\bf
    q}_{I,p,i+1}(T_p)={\bf 0},
  \label{equ:nonlinTrans}
\end{equation}
on the intervals $[T_{p},T_{p+1}]$, where the Jacobians
\begin{equation}
  {\bf J}_{i,p} = \frac{1}{T_{p+1}-T_p}\int_{T_p}^{T_{p+1}} \left
  .\frac{\partial {\bf N}({\bf q})}{\partial {\bf q}}\right|_{{\bf
      q}={\bf q}_{p,i}} \;\textrm{d}t.
  \label{equ:Jacobians}
\end{equation}
are now averaged on each time interval. Similarly, the homogeneous equations become
\begin{equation}
  \frac{\partial {\bf q}_{H,p,j,i+1}}{\partial t}=({\bf A}+{\bf
    J}_{i,p}){\bf q}_{H,p,j,i+1},
  \label{equ:nlHomIter}
\end{equation}
with the initial conditions
\begin{equation}
  {\bf q}_{H,p,j,i}(T_p)=
  \begin{cases}
    {\bf q}_{I,p-1,i+1}(T_p),  \;&j=0, p \neq 0,\\
    {\bf q}_{H,p-1,j-1,i+1}(T_p)  \;&1\leq j \leq p-1,\; p \neq 0,1.
  \end{cases}
\end{equation}
Once the homogeneous equations~(\ref{equ:nlHomIter}) are solved the
next iterate of the direct solution can be found via the sum
\begin{equation}
  {\bf q}_{p,i+1} = {\bf q}_{0}+{\bf q}_{I,p,i+1} + \sum_{j=1}^{p}
  {\bf q}_{H,p,j,+1},
  \label{equ:fullsolIter}
\end{equation}
on $[T_p,T_{p+1}]$. One way to determine if the full non-linear solution is converged, is by iterating until $\|{\bf q}_{i+1}-{\bf
  q}_{i}\|<\epsilon,$ where $\epsilon$ is a user-specified tolerance. We note that caution must be used to ensure that the solution is indeed converged or whether it has merely stagnated. Possible ways to help ensure proper convergence are to require that at least a certain amount of iterations are performed, or to randomly perturb the solution after some iterations to prevent the solution from getting trapped in a `plateau'. However, for this study we will not consider these technicalities. Once the non-linear solution
has been determined, the adjoint equation is then solved using the
algorithm presented in section~\ref{sec:linearAlg}.

\subsection{Governing equation solved in series}
\label{sec:hyb}

\begin{figure}
  \centering
  \begin{tabular}{cc}
    \begin{subfigure}[t]{0.4\textwidth}
      \includegraphics[width=\linewidth]{ 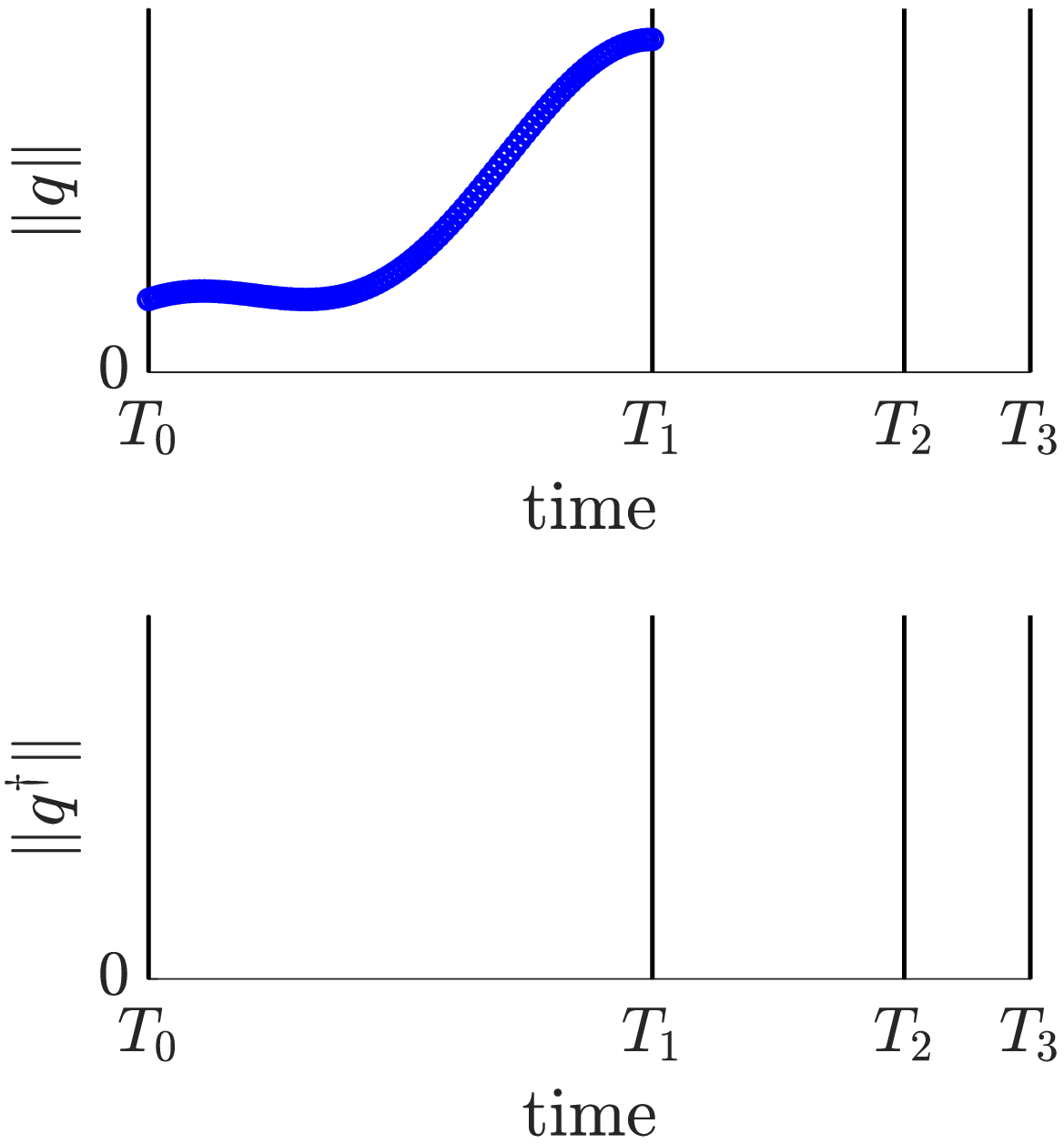}
      \caption{First, the direct equation is solved on the zeroth processor.\\~}
      \label{fig:hyb1}
    \end{subfigure}%
    & \hspace{1truecm}
    \begin{subfigure}[t]{0.4\textwidth}
      \includegraphics[width=\linewidth]{ 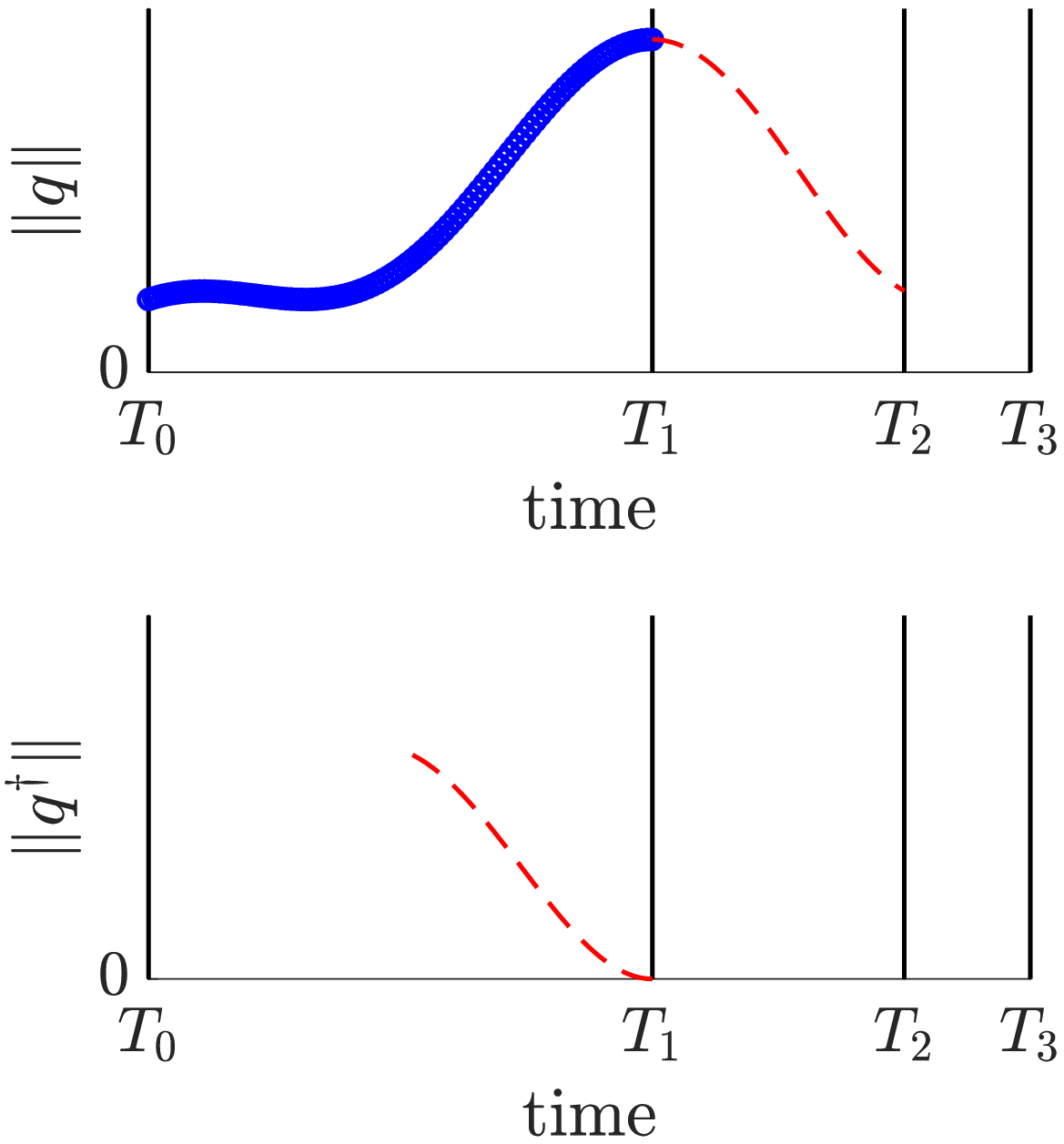}
      \caption{Whilst the first processor takes over the direct
        solution, the adjoint equation begins to be solved on the zeroth
        processor.}
  \label{fig:hyb2}
    \end{subfigure} \\
    \begin{subfigure}[t]{0.4\textwidth}
      \includegraphics[width=\linewidth]{ 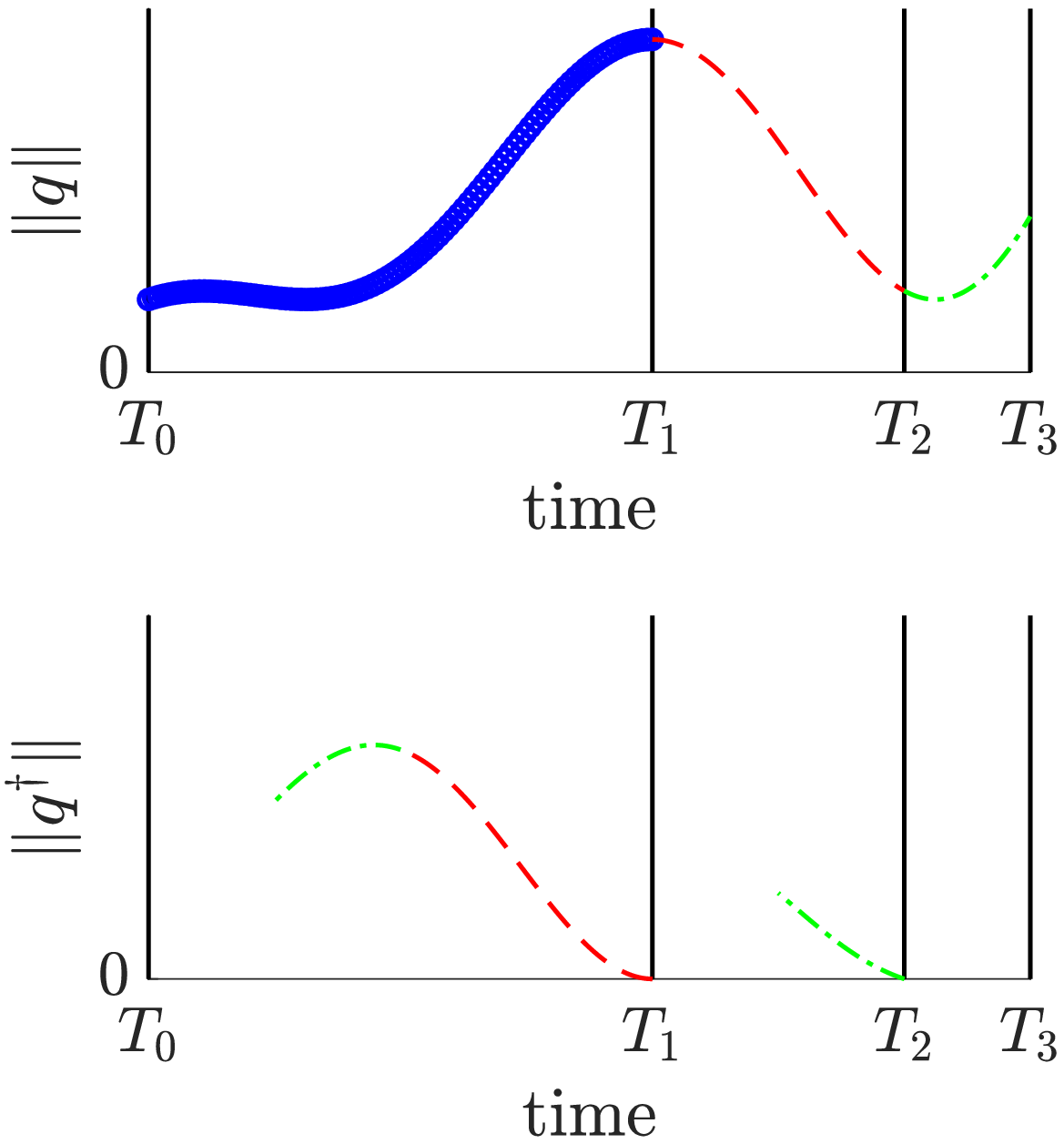}
      \caption{The direct solution is completed on the second
        processor. Simultaneously, the adjoint inhomogeneous solution
        is continued on processor zero, and begins on processor one.}
      \label{fig:hyb3}
    \end{subfigure}%
    & \hspace{1truecm}
    \begin{subfigure}[t]{0.4\textwidth}
      \includegraphics[width=\linewidth]{ 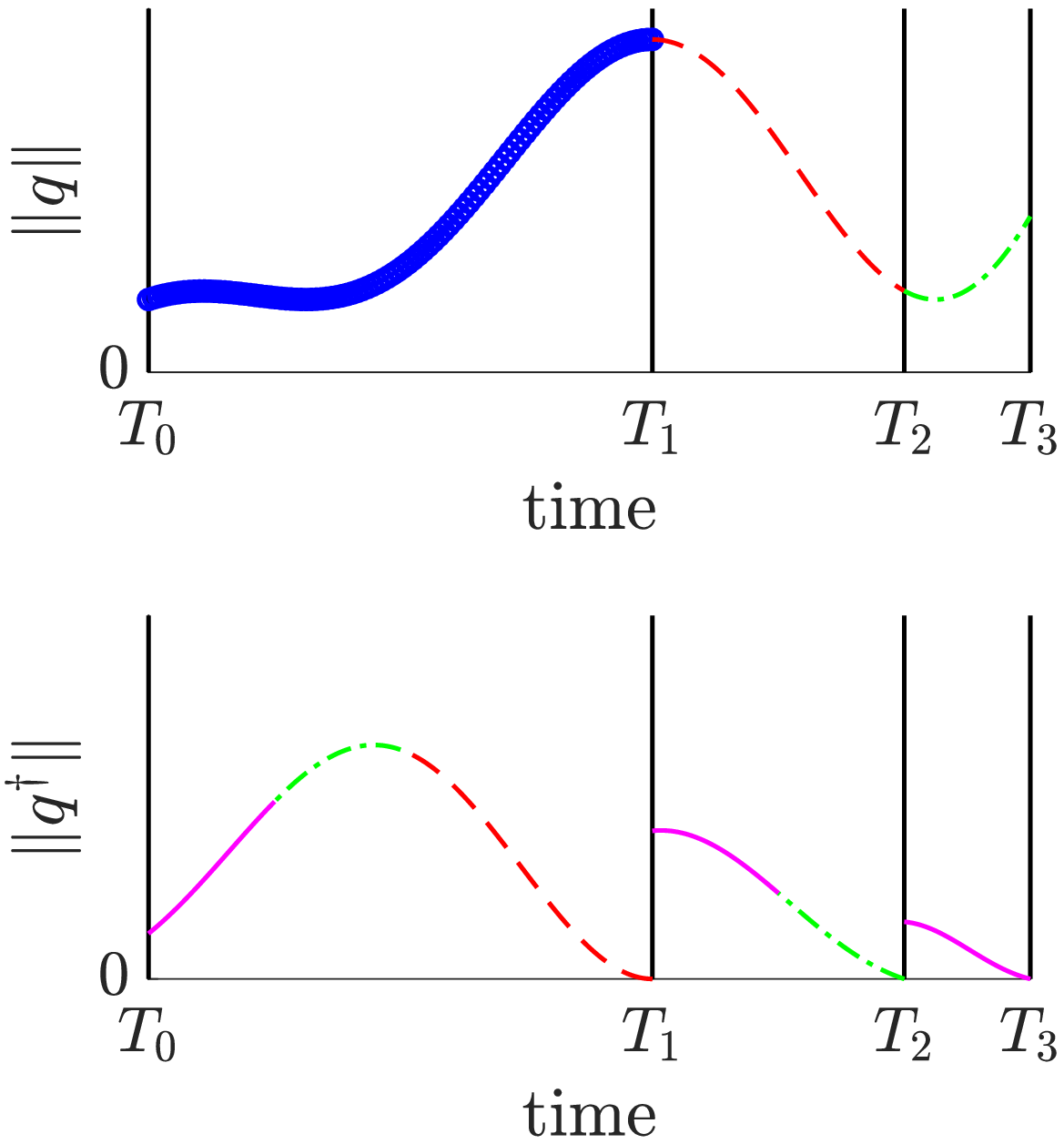}
      \caption{Finally, all adjoint inhomogeneous solutions are
        completed.\\~\\~}
      \label{fig:hyb4}
    \end{subfigure} 
  \end{tabular}
  \caption{The overlapping of the direct solution with the adjoint
    inhomogeneous solutions. The colours and different line styles
    show which parts of the solution are completed at the same time.}
  \label{fig:hybridGuide}
\end{figure}

Whilst the previous section shows that the governing equation can be
solved in parallel, problems may nevertheless arise. As there is no
guarantee that the iterations will converge, it is possible that the
direct solution may never be obtained. Additionally, even if the
iterations do converge, the number of iterations required might
outweigh the benefit obtained from solving the equation in
parallel. In these cases an alternative approach is needed; we propose
an algorithm in which the non-linear equation is solved in series,
whilst accommodating the time-parallelization of the adjoint
equation. This gives rise to the hybrid serial-direct-parallel-adjoint
method presented in this section and illustrated diagrammatically in figure \ref{fig:hybridGuide}.

For the direct equation we partition the time interval so that the
$p^\textrm{th}$ processor only solves equations on the interval
$[T_p,T_{p+1}]$. The direct equation is solved explicitly by
integrating the following equations in sequence
\begin{equation}
  \frac{\partial {\bf q}_p}{\partial t}=\boldsymbol{\mathcal{N}}({\bf q}_p,{\bf f}),
  \label{equ:seriesDirect}
\end{equation}
on $[T_p,T_{p+1}]$, with the initial conditions
\begin{equation}
  {\bf q}_p(T_p)=
  \begin{cases}
    {\bf q}_0,  \;& p = 0,\\
    {\bf q}_{p-1}(T_p) \;&1\leq p \leq N-1,
  \end{cases}
\end{equation}
where the notation ${\bf q}_p$ is used to denote the full solution
${\bf q}$ solved on processor $p$. As the $p^\textrm{th}$ processor
needs an initial condition from the $(p-1)^\textrm{th}$ processor,
later processors must wait until previous processors have finished
their direct solutions. However, unlike in the previous algorithms,
the full direct solution is available on the interval $[T_p,T_{p+1}]$
as each processor $p$ finishes solving
equation~(\ref{equ:seriesDirect}). Therefore, while later processors
are computing the direct solutions, earlier processors can begin
solving the adjoint inhomogeneous equations. The
algorithm proceeds in this manner with the direct solution being
solved across all the processors, with the adjoint inhomogeneous
equations beginning on each processor once the
direct solution has been computed. After the adjoint inhomogeneous
equations have been solved, all that remains is to compute the adjoint
homogeneous equations.

In the algorithms presented previously, the initial condition for the adjoint variable was absorbed into the adjoint inhomogeneous equation. This was possible due to the adjoint solution occurring after the direct solution was calculated, and therefore the adjoint final time condition being available. However, in overlapping the direct and adjoint inhomogeneous solutions, this is not possible as the adjoint inhomogeneous solutions occur before this final time condition can be calculated. Therefore, in this case the inhomogeneous equations that must be solved are
\begin{equation}
  \frac{\partial {\bf q}^\dagger_{I,p}}{\partial t}=-{\bf A}^H{\bf
    q}^\dagger_{I,p}-\overline{\frac{\partial \mathcal{J}_I}{\partial
      {\bf q}}}, \qquad {\bf q}^\dagger_{I,p}(T_{p+1})={\bf 0},
      \label{equ:hybadjinhom}
\end{equation}
on $[T_{p+1},T_{p}]$, i.e. the adjoint inhomogeneous equations~(\ref{equ:adjInhom}) without absorbing the final time condition into them. The fact that we cannot include the effect of an adjoint final time condition directly in the equations means that we have to solve one more homogeneous problem stemming from this final time condition. Hence, in addition to solving the adjoint homogeneous problems~(\ref{equ:adjHomNormal}), we must additionally solve
\begin{equation}
  \frac{\partial {\bf q}^\dagger_{H,N-1}}{\partial t}=-{\bf A}^H{\bf
    q}^\dagger_{H,N-1}, \qquad {\bf q}^\dagger_{H,N-1}(T)={\bf q}^\dagger(T),
    \label{equ:adjhomfinal}
\end{equation}
on $[T,0]$ if there is a non-zero final time condition.

In contrast to the previous algorithms in which an equidistant time
partition would be a sensible choice, the hybrid algorithm requires
slightly more attention. Using the hybrid scheme means that the
inhomogeneous adjoint equations start on earlier processors first. If
an equidistant time partition were used, then these earlier processors
would terminate first and consequently would idle for information from
subsequent processors to begin solving the homogeneous adjoint
equation. This time spent waiting limits the efficiency of the
algorithm, as unused computing power is squandered. Therefore, it is
imperative to choose the time partitioning strategically, thus
avoiding the inefficiency of the equally-spaced time grid.

To prevent earlier processors waiting on later processors, all the
inhomogeneous adjoint equations being solved need to terminate at the
same time, i.e., the time for the full solution to be completed on one
processor after having completed on the previous ones plus the
subsequent adjoint inhomogeneous evaluation should remain constant on
each processor. By introducing the time $\tau_I$ taken per time unit
for the direct solve and similarly the time $\tau_I^\dagger$ taken per
time unit for the adjoint inhomogeneous solve, we can write this
relationship in terms of the time partition for processors $n$ and
$n+1$ as
\begin{eqnarray}
  \underbrace{\tau_I\left(\sum_{i=0}^n(T_{i+1}-T_i)\right)}_{\text{direct}}
  + \underbrace{\tau_I^\dagger(T_{n+1}-T_n)}_{\text{adjoint}} &=&
  \underbrace{\tau_I\left(\sum_{i=0}^{n+1}(T_{i+1}-T_i)\right)}_{\text{direct}}\nonumber\\
  &&+ \underbrace{\tau_I^\dagger(T_{n+2}-T_{n+1})}_{\text{adjoint}}.
\end{eqnarray}
This relationship can be rearranged to produce
\begin{equation}
  T_{n+2}-T_{n+1}=\frac{k}{k+1}(T_{n+1}-T_{n}).
  \label{equ:paritionRatio}
\end{equation}
 with $k = \tau_I^\dagger/\tau_I$. This expression necessarily implies
 that the simultaneous termination of the inhomogeneous evaluations is
 achieved by making subsequent time partitions finer as we approach
 the final time. Equation~(\ref{equ:paritionRatio}) is a second-order
 difference equation, with boundary conditions $T_0=0$, and $T_{N}=T$
 which can be solved to yield
\begin{equation}
  T_n =
  T\left(\frac{1-\left(\displaystyle{\frac{k}{k+1}}\right)^{n}}
  {1-\left(\displaystyle{\frac{k}{k+1}}\right)^{N}}\right),
  \label{equ:nonEquiP}
\end{equation}
giving an analytic expression for the time partition in terms of the
ratio $k$, which can be measured numerically.

One important difference of this hybrid serial-direct-parallel-adjoint
algorithm to the previous approaches must be pointed out. In the
previous cases, an equidistant time partition can be used. This means
that, assuming the adjoint homogeneous
equations~(\ref{equ:adjHomNormal}) take the same amount of time to be
solved on each processor, that they can be solved by instead solving
the series of equations~(\ref{equ:adjHomLocal}). This keeps the solved
equations local to each processor's time partition, meaning that the
direct solution does not need to be distributed to other
processors. However, in this serial-direct-parallel-adjoint approach
the time partition is inherently non-equidistant. Solving the adjoint
homogeneous part using the series of equations~(\ref{equ:adjHomLocal})
will mean that later processors, which have a smaller time partition,
will need to wait to pass their states to lower processors which solve
equations on increasingly longer time intervals. The thus introduced
lag scales with the number of processors, negating any speedup
obtained from overlapping the direct and adjoint solutions. In order
to circumvent this problem, the adjoint homogeneous equations must
instead be solved in their original form given by
equations~(\ref{equ:adjHomNormal}) and~(\ref{equ:adjhomfinal}). Although this means that the direct
solution must be distributed to all processors, this is a one-time
cost that must be paid, as it enables the algorithm to achieve an
overall speedup. We choose to solve the $j^\textrm{th}$ adjoint
homogeneous equation on processor $p$. The
full adjoint solution on $[T_p,T_{p+1}]$ can then be found by the sum
\begin{equation}
  {\bf q}^\dagger = {\bf q}^\dagger(T)+ {\bf q}_{I,p}^\dagger +
  \sum_{j=0}^{N-1} {\bf q}_{H,j}^\dagger.
  \label{equ:fullAdjSolHyb}
\end{equation}

Figure~\ref{fig:hybridGuide} illustrates how the hybrid algorithm
plays out for a three processor arrangement. The algorithm begins by
solving the direct equation for processor zero
(figure~\ref{fig:hyb1}). The integration of the direct equation then
continues on the next processor in figure~\ref{fig:hyb2} whilst the
inhomogeneous adjoint equation begins for $p=0$. Similarly,
figure~\ref{fig:hyb3} shows that, when the direct integration finishes
for $p=1,$ the direct integration continues on processor two, whilst
the inhomogeneous integration begins on processor one and continues on
processor zero. Finally, figure~\ref{fig:hyb4} shows that after the
direct integration has finished, all inhomogeneous adjoint equations
terminate at the same time. The adjoint solution can then be completed
by solving the homogeneous equations~(\ref{equ:adjHomNormal}) and (\ref{equ:adjhomfinal}) in
parallel (not shown).

\subsection{Step-by-step summary of the algorithms and scaling analysis}

Now that the three algorithms have been presented, we provide a
step-by-step summary for ease of implementation. A scaling analysis is
also conducted to assess the theoretical speedups possible.

\subsubsection{Linear governing equation}
\label{sec:linEquSteps}

\begin{algorithm}
  \caption{Direct-adjoint loop solved in parallel for a linear governing equation}
  \label{alg:lin}
  \begin{algorithmic}[1]
    \Procedure{DirectAdjointLoop}{${\bf q}_0$,N} \Comment{Perform a single direct-adjoint loop in parallel with $N$ processors}
    \State Partition the time domain $[0,T]$ into $N$ equal parts $0=T_0<T_1<...<T_N=T$
    \For{\textbf{processor} $p \in \{0,...N-1\}$}
    \State Solve equation (\ref{equ:linInhomEqus}) on $[T_p,T_{p+1}]$ \Comment{Direct inhomogeneous equations}\label{alg:LinInhom}
    \If{$p<N-1$}
    \State Send the state ${\bf q}_{I,p}(T_{p+1})$ to processor $p+1$\Comment{communication}
    \EndIf
    \If{$p>0$}
    \State Receive the state ${\bf q}_{I,p-1}(T_{p})$ from processor $p-1$\Comment{communication}
    \EndIf
    \For{$j \in \{0,...,p-1\}$} \label{alg:Linhom}
    \If{$j>0$}
    \State Receive the state ${\bf q}_{H,p-1,j-1}(T_{p})$ from processor $p-1$
    \EndIf
    \State Solve equation (\ref{equ:adjHomLocal}) on $[T_p,T_{p+1}]$ \Comment{Direct homogeneous equations}
    \If{$j<p-1$}
    \State Send the direct homogeneous state ${\bf q}_{H,p,j}(T_{p+1})$ to processor $p+1$ \Comment{communication}
    \EndIf
    \EndFor
    \State Form the full direct solution on $[T_p,T_{p+1}]$ by the sum (\ref{equ:fullsol0})
    \State Initialize the adjoint state with the final time condition (\ref{equ:adjFinalTime})
    \State Solve the adjoint inhomogeneous equation  (\ref{equ:adjInhom}) on $[T_{p+1},T_p]$ \Comment{Adjoint inhomogeneous equations}
    \If{$p>0$}
    \State Send the state ${\bf q}^\dagger_{I,p}(T_{p})$ to processor $p-1$ \Comment{communication}
    \EndIf
    \If{$p<N-1$}
    \State Receive the state ${\bf q}^\dagger_{I,p+1}(T_{p+1})$ from processor $p+1$ \Comment{communication}
    \EndIf
    \algstore{linalg}
\end{algorithmic}
\end{algorithm}
    
\begin{algorithm}                     
\begin{algorithmic}[1]                   
\algrestore{linalg}
    \For{$j \in \{0,...,N-1-p\}$}
    \If{$j>0$}
    \State Receive the state ${\bf q}^\dagger_{H,p+1,j-1}(T_{p+1})$ from processor $p+1$
    \EndIf
    \State Solve equation (\ref{equ:adjHomLocal}) on $[T_{p+1},T_{p}]$ \Comment{Adjoint homogeneous equations}
    \If{$j<N-1-p$}
    \State Send the state ${\bf q}^\dagger_{H,p,j}(T_p)$ to processor $p-1$ \Comment{communication}
    \EndIf
    \EndFor
    \State Form the full adjoint solution on $[T_p,T_{p+1}]$ by the sum (\ref{equ:fulladjsol})
    \EndFor 
    \EndProcedure
  \end{algorithmic}
\end{algorithm}

Algorithm~\ref{alg:lin} summarizes the linear case discussed in
section~\ref{sec:linearAlg}. To quantify the speedup for the
algorithm, we will introduce the time taken per time unit to solve the
direct inhomogeneous and adjoint equations, denoted $\tau_I$ and
$\tau_H$, respectively. Similarly, we introduce adjoint counterparts
$\tau^\dagger_I$, $\tau^\dagger_H$. For every processor, we assume
that step~\ref{alg:LinInhom} of algorithm~\ref{alg:lin} takes the same
amount of time to solve and occurs perfectly in sync. This means that
the time taken to solve the inhomogeneous direct equations is
$(T/N)\tau_I$. After the direct inhomogeneous equations are solved,
the direct homogeneous equations must be solved. This is achieved by
performing the for-loop beginning in step~\ref{alg:Linhom}. Although
each processor solves a different amount of direct homogeneous
integrations, it is processor $N-1$ that solves the most. Therefore,
the time for this loop is $(T(N-1)/N)\tau_H$. Hence, the total time to
solve the direct equation in parallel is
\begin{equation}
  T_D = \frac{T}{N}\tau_I+\frac{T(N-1)}{N}\tau_H.
  \label{eq:timeDirectPara}
\end{equation}
Performing the same analysis on the adjoint part of the algorithm
gives the time taken to solve the adjoint equation as
\begin{equation}
  T_A = \frac{T}{N}\tau_I^\dagger+\frac{T(N-1)}{N}\tau_H^\dagger.
  \label{equ:adjointTime}
\end{equation}
Adding the two parts together gives the total time for one
direct-adjoint loop $T_L=T_D+T_A$ as
\begin{equation}
  T_L = \frac{T}{N}(\tau_I+\tau_I^\dagger) +
  \frac{T(N-1)}{N}(\tau_H+\tau_H^\dagger).
\end{equation}
In series, the time taken for the algorithm is simply
$T_L^S=T(\tau_I+\tau_I^\dagger)$. Hence, the speedup $s=T^S_L/T_L$
(shown as the reciprocal speedup for convenience) is
\begin{equation}
  s^{-1} = \frac{1}{N} +
  \frac{(N-1)}{N}\left(\frac{\tau_H+\tau_H^\dagger}
       {\tau_I+\tau_I^\dagger}\right).
       \label{equ:LinSpeedup}
\end{equation}
The speedup~(\ref{equ:LinSpeedup}) shows us that there is a
theoretical maximum to the possible speedup. Indeed, as
$N\rightarrow\infty$, $s\rightarrow
(\tau_I+\tau_I^\dagger)/(\tau_H+\tau_H^\dagger)$. This highlights the
importance of solving the homogeneous equation faster than the
inhomogeneous equation. The faster the homogeneous integration over
the inhomogeneous evaluation, the more efficient the parallel-in-time
approach becomes.

\subsubsection{Non-linear governing equation solved in parallel}

\begin{algorithm}
  \caption{Direct-adjoint loop solved in parallel for a non-linear governing equation}
  \label{alg:nonlin}
  \begin{algorithmic}[1] 
    \Procedure{DirectAdjointLoop}{${\bf q}_0$,N} \Comment{Perform a single direct-adjoint loop in parallel with $N$ processors}
    \State Partition the time domain $[0,T]$ into $N$ equal parts $0=T_0<T_1<...<T_N=T$
    \State Choose a tolerance \textit{tol}
    \State $\textit{err}=1e^{10}$
    \State $\textit{iter}=0$
    \For{\textbf{processor} $p \in \{0,...N-1\}$}
    \State Make a guess for the direct solution ${\bf q}_{p,0}$
    \State Form the Jacobian terms by equation (\ref{equ:Jacobians})
    \While{$\textit{err}>\textit{tol}$}\label{alg:nonLinWhile}
    \State Solve equation (\ref{equ:nonlinTrans}) on $[T_p,T_{p+1}]$ \Comment{Direct inhomogeneous equations}
    \If{$p<N-1$}
    \State Send the state ${\bf q}_{I,p,\textit{iter}+1}(T_{p+1})$ to processor $p+1$ \Comment{communication}
    \EndIf
    \If{$p>0$}
    \State Receive the state ${\bf q}_{I,p-1,\textit{iter}+1}(T_{p})$ from processor $p-1$ \Comment{communication}
    \EndIf
    \For{$j \in \{0,...,p-1\}$}
    \If{$j>0$}
    \State Receive the state ${\bf q}_{H,p-1,j-1,\textit{iter}+1}(T_{p})$ from processor $p-1$ \Comment{communication}
    \EndIf
    \State Solve equation (\ref{equ:nlHomIter}) on $[T_p,T_{p+1}]$ \Comment{Direct homogeneous equations}
    \If{$j<p-1$}
    \State Send the state ${\bf q}_{H,p,j,\textit{iter}+1}(T_{p+1})$ to processor $p+1$ \Comment{communication}
    \EndIf
    \EndFor
    \State Form the current iterate of the full direct solution ${\bf q}_{p,\textit{iter}+1}$ on $[T_p,T_{p+1}]$ by the sum (\ref{equ:fullsolIter})
    \State Form the Jacobian terms by equation (\ref{equ:Jacobians}) \label{alg:nonLinJac}
    \State Calculate the error $\textit{err}$
    \State $\textit{iter}\leftarrow \textit{iter}+1$
    \EndWhile
    \algstore{nlalg}
\end{algorithmic}
\end{algorithm}

\begin{algorithm}                     
\begin{algorithmic} [1]                   
\algrestore{nlalg}
    \State Initialize the adjoint state with the final time condition (\ref{equ:adjFinalTime})
    \State Solve equation (\ref{equ:adjInhom}) on $[T_{p+1},T_p]$ \Comment{Adjoint inhomogeneous equations}
    \If{$p>0$}
    \State Send the state ${\bf q}^\dagger_{I,p}(T_{p})$ to processor $p-1$ \Comment{communication}
    \EndIf
    \If{$p<N-1$}
    \State Receive the state ${\bf q}^\dagger_{I,p+1}(T_{p+1})$ from processor $p+1$ \Comment{communication}
    \EndIf
    \For{$j \in \{0,...,N-1-p\}$}
    \If{$j>0$}
    \State Receive the state ${\bf q}^\dagger_{H,p+1,j-1}(T_{p+1})$ from processor $p+1$ \Comment{communication}
    \EndIf
    \State Solve equation (\ref{equ:adjHomLocal}) on $[T_{p+1},T_{p}]$ \Comment{Adjoint homogeneous equations}
    \If{$j<N-1-p$}
    \State Send the state ${\bf q}^\dagger_{H,p,j}(T_p)$ to processor $p-1$ \Comment{communication}
    \EndIf
    \EndFor
    \State Form the full adjoint solution on $[T_p,T_{p+1}]$ by the sum (\ref{equ:fulladjsol})
    \EndFor 
    \EndProcedure
  \end{algorithmic}
\end{algorithm}

The case of a non-linear governing equation, see
section~\ref{sec:nonLinParallel}, solved in parallel is summarized
in algorithm~\ref{alg:nonlin}. It has a strong resemblance to
algorithm~\ref{alg:lin} except for the fact that the direct equation
is now solved in an iterative fashion with the while-loop starting on
line~\ref{alg:nonLinWhile}. For the purposes of assessing the speedup,
we will include the cost of calculating the Jacobian terms on
line~\ref{alg:nonLinJac} with the direct inhomogeneous
equation. Furthermore, as the inhomogeneous direct equation in
parallel is not the same as the direct equation in series, we
introduce the time taken per time unit $\tau_D^S$ for the direct
equation in series. This gives the time taken for a direct-adjoint
loop in series as $T_S=T(\tau_D^S+\tau_I^\dagger)$.

We introduce here the number of iterations $K(N)$ we require to
achieve convergence for our non-linear scheme. Note the dependence of
the iterations on the number of processors; the fineness of the time
partition will affect the number of iterations. As each iteration
involves the same steps as the linear case we obtain that the time
taken for the direct equation is simply $K(N)$ times
equation~(\ref{eq:timeDirectPara}), i.e.,
\begin{equation}
  T_D = K(N)\left(\frac{T}{N}\tau_I + \frac{T(N-1)}{N}\tau_H\right).
\end{equation}
The time taken to solve the adjoint equation is still given by
equation~(\ref{equ:adjointTime}), due to the required steps being
identical to algorithm~\ref{alg:lin}. Hence, the time taken to perform
one direct-adjoint loop in parallel is
\begin{equation}
    T_L=K\left(\frac{T}{N}\tau_I + \frac{T(N-1)}{N}\tau_H\right) +
    \frac{T}{N}\tau_I^\dagger + \frac{T(N-1)}{N}\tau_H^\dagger,
\end{equation}
giving the reciprocal speedup as 
\begin{equation}
  s^{-1} = \frac{1}{N}\left(\frac{K\tau_I+\tau_I^\dagger}
  {\tau_D^S+\tau_I^\dagger}\right) +
  \frac{(N-1)}{N}\left(\frac{K\tau_H+\tau_H^\dagger}
       {\tau_D^S+\tau_I^\dagger}\right).
  \label{equ:nonLinSpeedup}
\end{equation}

Equivalently to the linear case, equation~(\ref{equ:nonLinSpeedup})
shows the importance of the homogeneous solvers being faster than the
inhomogeneous ones. However, unlike the linear case, the factor of $K$
in the speedup equation implies that small $N$ can lead to $s<1$,
signifying that no speedup is possible using a parallel
approach. Requiring the speedup to be greater than one produces the
condition that
\begin{equation}
  N > \frac{K(\tau_I-\tau_H)+\tau_I^\dagger-\tau_H^\dagger}
  {\tau_D^S+\tau_I^\dagger-\tau_H-K\tau_H^\dagger}.
  \label{equ:speedupCond}
\end{equation}
Applying condition~(\ref{equ:speedupCond}) to the linear case, where
$K=1$ and $\tau_D^S=\tau_I$ simply gives the condition $N>1$, implying
that a speedup is achievable for any number of processors. In the
non-linear case, this no longer holds true, and $N$ must
satisfy~(\ref{equ:speedupCond}) in order to be faster than the serial
case. Once again, the algorithm has a theoretical maximum speedup; as
$N\rightarrow\infty$,
$s\rightarrow(\tau_D^S+\tau_I^\dagger)/(K\tau_H+\tau_H^\dagger)$. This
is similar to the linear case, except for the presence of a factor $K$
in the denominator which implies that as the number of iterations
increases, the maximum speedup decreases.

\subsubsection{Governing equation solved in series}

\begin{algorithm}
  \caption{Direct-adjoint loop solved in serial-direct-parallel-adjoint for a non-linear governing equation}
  \label{alg:hyb}
  \begin{algorithmic}[1] 
    \Procedure{DirectAdjointLoop}{${\bf q}_0$,N} \Comment{Perform a single direct-adjoint loop in parallel with $N$ processors}
    \State Estimate the value of $k$ in (\ref{equ:paritionRatio}) by performing short direct and adjoint homogeneous calculations.
    \State Partition the time domain $[0,T]$ into $N$ parts $0=T_0<T_1<...<T_N=T$ using equation (\ref{equ:nonEquiP})
    \For{\textbf{processor} $p \in \{0,...N-1\}$}
    \If{$p>0$}
    \State Wait to receive the state ${\bf q}_{p-1}(T_{p})$ from processor $p-1$ \Comment{communication}
    \EndIf
    \State Solve equation (\ref{equ:nonlinTrans}) on $[T_p,T_{p+1}]$ \Comment{Direct equation}
    \If{$p<N-1$}
    \State Send the direct state ${\bf q}_{p}(T_{p+1})$ to processor $p+1$ \Comment{communication}
    \EndIf
    \State Solve the adjoint inhomogeneous equation (\ref{equ:hybadjinhom}) on $[T_{p+1},T_p]$ \Comment{Adjoint inhomogeneous equations}
    \EndFor
    \State Scatter the direct solution to all processors \Comment{communication}
    \For{\textbf{processor} $p \in \{0,...N-1\}$}
    \If{$p<N-1$}
    \State Solve the $j=p^\text{th}$ adjoint homogeneous equation
    (\ref{equ:adjHomNormal}) on $[T_{p+1},0]$ \Comment{Adjoint homogeneous equations}
    \EndIf
    \If{$p=N-1$ {\bf and} $\mathcal{J}_T\neq 0$  }
    \State Solve the $j=N-1^\text{th}$ adjoint homogeneous equation
    (\ref{equ:adjhomfinal}) on $[T,0]$ \Comment{Adjoint homogeneous equations}
    \EndIf
    \EndFor
    \State Form the full adjoint solution with the sum (\ref{equ:fullAdjSolHyb})
    \Comment{communication}
    \EndProcedure
  \end{algorithmic}
\end{algorithm}
Lastly, we consider the hybrid direct-serial-parallel-adjoint case
introduced in section~\ref{sec:hyb}. The procedure is summarized in
algorithm~\ref{alg:hyb}. To obtain the speedup, we initially consider
the time needed to solve the direct equation together with the adjoint
inhomogeneous equations. Due to the overlapping nature of the
algorithm, and the way in which we choose the non-equidistant
partitioning, all inhomogeneous parts of the adjoint equation should
finish at the same time. First, the direct equation is integrated in
series giving a time of $\tau_IT$. Following this, the adjoint
inhomogeneous equations must all be integrated for $(T_N-T_{N-1})$
time units. The time taken to finish the direct equation and the
inhomogeneous adjoint equations is then
$\tau_IT+(T_N-T_{N-1})\tau_I^\dagger$. Finally, the time taken to
finish the direct-adjoint loop is equivalent to the time taken to
solve the longest homogeneous equation, which is solved for either $T_{N-1}$ or $T$
time units depending on whether there is an adjoint final time condition. By denoting this homogeneous time as $T_H$ we obtain a total time for the direct-adjoint loop as
\begin{equation}
  T_L = T\tau_I+(T_N-T_{N-1})\tau_I^\dagger+T_H\tau_H^\dagger.
\end{equation}
Dividing this expression by the time it takes to solve the system in
series ($T_L^S$ from the linear discussion in
section~(\ref{sec:linEquSteps})) gives the reciprocal speedup as
\begin{equation}
  s^{-1}=\frac{T\tau_I+(T_N-T_{N-1})\tau_I^\dagger +
    T_H\tau_H^\dagger}{T(\tau_I+\tau_I^\dagger)}.
  \label{equ:hybridSpeedup}
\end{equation}

Whilst this theoretical speedup is considerably more complicated than
the other cases, we can still consider the maximum speedup
achievable. We note that with the time partitioning based on
equation~(\ref{equ:nonEquiP}), we obtain $T_{N-1}\rightarrow T$ and
$(T_{N}-T_{N-1})\rightarrow 0$ as $N\rightarrow \infty$. As $T_H$ is either $T_{N-1}$ or $T$ we obtain $s\rightarrow (\tau_I+\tau_I^\dagger)/(\tau_I+\tau_H^\dagger)$,
showing that effectively we achieve a speedup by replacing an
inhomogeneous adjoint equation with a homogeneous one.

\subsection{Checkpointing}\label{sec:checkpointing}

At this point it is important to note that, owing to the addition of
the forcing term in the adjoint formulation, there is an explicit
dependence of the adjoint equations on the direct
variables. Therefore, the direct variables need to be stored at all
time steps during the forward sweep and injected, at the appropriate
time, into the adjoint equations. For high-resolution cases and large
time horizons, we cannot afford to store all necessary direct
variables, as the memory requirements may exceed the necessary amount
of RAM of our computational system. In this case a checkpointing
scheme is needed.

\begin{figure}
  \center
  \includegraphics[width=0.7\textwidth]{ 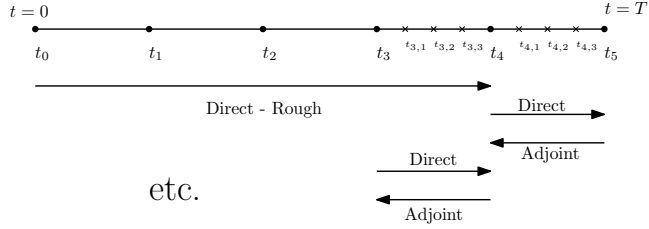}
  \caption{Graphical representation of the checkpointing procedure.}
  \label{Fig:CheckPointing}
\end{figure}

A graphical representation of the checkpointing process, where four
checkpoints $t_1,...,t_4$ are used, can be seen in
figure~\ref{Fig:CheckPointing}. Initially, the direct equation is
integrated to the final checkpoint $t_4$ without saving any
intermediate states (this run is denoted \textit{Rough} in the
figure). The direct equation is then integrated from the final
checkpoint to the final time saving all intermediate states. As the
full solution is now available between $t_4$ and $t_5,$ the adjoint
variable can be initialized with the relevant final time condition and
integrated back to $t_4$. The intermediate states between $t_4$ and
$t_5$ are now deleted, freeing up the RAM. This allows the direct
equation to be integrated from the penultimate checkpoint $t_3$ to
$t_4,$ again with all intermediate states saved. The adjoint equation
can then continue its backwards integration from $t_4$ to $t_3$. By
continuing this checkpointing scheme backwards in this manner, the
adjoint equation can be integrated back to the initial time at a
fraction of the memory cost.

In view of our parallel-in-time algorithms we incorporate a
checkpointing scheme as follows. We begin by defining the chosen
location of our checkpoints. Then, the direct or direct-adjoint loops
between each checkpoints are carried out in parallel. In this way the choice of checkpoints is independent of the choice of time-partition. For example, in
the four-checkpoint example we would first carry out the
direct-\textit{Rough} integration by integrating in parallel (or
serially in the hybrid case) the equation in the interval $[t_0,t_1]$,
followed by the interval $[t_1,t_2]$ and so on, until the
direct-\textit{Rough} integration is finished. When this is
accomplished, the direct-adjoint loops are performed using the
parallel algorithms developed so far on the checkpointed intervals,
preceding to the first time interval. Although the checkpoint times
can be chosen optimally~\citep[see][]{revolve}, for our numerical implementation we proceed by using
equispaced checkpoints for ease of implementation.

\section{Implementation details}

As revealed by the scaling analysis, it is critically important that
the homogeneous equations are solved faster than the inhomogeneous
ones. Indeed, we can see from the theoretical
scalings~(\ref{equ:LinSpeedup}), (\ref{equ:nonLinSpeedup})
and~(\ref{equ:hybridSpeedup}) that the maximum speedup possible for
each algorithm crucially depends on solving the homogeneous equation
more efficiently. This implies a natural restriction on the kind of
equations that we are able to parallelize in time with this
approach. Utilizing the approach of~\citet{GanderParaExp}, we identify
stiff equations as good candidates for time parallelization. The
homogeneous components that arise from these stiff problems can be
solved efficiently and quickly with exponential time-steppers, an
example of which was derived and applied to stiff equations
by~\citet{cox2002}. As exponential time-steppers are based on the
exact integration of the homogeneous equation, via the matrix
exponential or approximations thereof, they alleviate the stiffness of
the problem and consequently are able to take larger time-steps than
traditional time-stepping methods based on linear multistep or
multistage formula. This section begins with a brief discussion of
exponential integrators. Next, we demonstrate the application of each
of the three algorithms developed to a direct-adjoint looping
procedure. Special attention is given to how well the scaling of the
algorithm agrees with the ones derived theoretically, and how the
maximum speedup depends on the stiffness of the problem.

\subsection{Exponential integrators}

Exponential integration is based on the analytic solution to a linear,
autonomous, homogeneous problem in terms of the matrix
exponential. For the linear problem $\dot{\bf y}={\bf A}{\bf y}$ with
initial condition ${\bf y}(0)={\bf y}_0$, the solution at time $\Delta
t$ is ${\bf y}=\exp(\Delta t {\bf A}){\bf y}_0$. Exponential
integrators rely on approximating the matrix-vector product
$\exp(\Delta t {\bf A}){\bf v}$, which advances the state ${\bf v}$ by
$\Delta t$ time units. Even though this approach is based on the exact
integration for a linear equation, this approach can still be used for
non-linear homogeneous systems, albeit with smaller time-steps. In
either case, exponential time-steppers can perform better than
alternative methods for stiff problems~\citep[see, e.g.,][]{kassam2005}. A short overview of two possible exponential integrators is now given. For a more thorough review, see the recent review paper by~\citet{matrixFunctions}.

One way to approximate the matrix-vector product $\exp(\Delta t
{\bf A}){\bf v}=\exp(
{\bf B}){\bf v}$ is via the orthogonalized, order-$m$ Krylov subspace, i.e. by forming the subspace
\begin{equation}
    \mathcal{K}_m({\bf B},{\bf v})=\{{\bf v},{\bf B}{\bf v},...,{\bf B}^{m-1}{\bf v}\},
\end{equation}
and orthogonalizing it~\citep{krylovApprox}. By computing this subspace we can obtain a lower-dimensional representation
of the matrix, namely ${\bf B}{\bf v}\approx{\bf V}_m{\bf H}_m{\bf V}_m^T{\bf v}$
where ${\bf H}_m$ is an $m \times m$ matrix. The product of the matrix exponential and a vector can then be approximated as  $\exp({\bf
  B}){\bf v}\approx{\bf V}_m\exp({\bf H}_m){\bf V}_m^T{\bf v}$, which can be computed
significantly more quickly and efficiently if $m$ is far smaller
than the dimension of ${\bf B}$. If the Krylov
representation is accurate for small $m$, then this method provides an
efficient way of calculating the matrix-exponential-vector product. However, if $m$
needs to be large (e.g., if ${\bf B}$ has a large spectral radius) then
this method becomes time-consuming and memory-intensive, as $m$
Krylov vectors must be calculated and stored.

An alternative to the Krylov-based method is based on Newton
interpolation of the exponential function
$\phi(z)=(\exp(z)-1)/z$~\citep{caliari2007}. Once the interpolation is
found, the matrix exponential $\exp({\bf B})$ can easily be retrieved
from $\phi({\bf B})$. When computing the Newton interpolation the
choice of interpolation points $\{\zeta_j\}_{j=0}^N$ is critically
important, as a clever choice of $\zeta_j$ reduces the amount of
points ($N+1$) needed for an accurate interpolation. One way to
achieve superlinear convergence of the matrix polynomial to the matrix
exponential ${\phi}({\bf Z})$ is to use Chebyshev nodes on the real
focal interval $[a,b]$, computed such that `the ``minimal'' ellipse of
the confocal family that contains the spectrum (or the field of
values) of the matrix is not too ``large''
'~\citep{realLeja}. However, the locations of Chebyshev nodes are
dependent on the interpolation degree $N+1$. Therefore, if the
interpolation degree needs to be increased in order to achieve a
user-specified error, then the entire interpolation must be
recomputed. This cost can be circumvented by instead using real Leja
points, which still achieve a superlinear convergence. The advantage
of Leja points over Chebyshev points is that they do not depend on the
degree of interpolation. Hence, if a smaller error is required, the
degree can be increased without recomputing the interpolation. Newton
interpolation of the exponential function $\phi(z)$ is called the real
Leja points method (ReLPM) and is summarized by~\citet{realLeja}. In
the case of large sparse matrices the ReLPM can become more
advantageous over Krylov-based methods due to its lower memory
cost. We choose to use the ReLPM method for our implementation so that
the developed code is more easily portable to future large-scale
problems. For a more detailed discussion of exponential integrators
see~\citet{skene}, which also contains a preliminary version of this
work.

\subsection{Linear governing equation}

To test the linear algorithm~\ref{alg:lin}, we consider a forced
advection-diffusion equation.

\begin{equation}
  \frac{\partial q}{\partial t}=-a\frac{\partial q}{\partial
    x}-a\frac{\partial q}{\partial y} +D\nabla^2 q+ f\sin(\omega t),
\end{equation}
where the speed of advection is governed by $a$ and the diffusion
coefficient is $D$. We take the domain to be
$\Omega=[0,2\pi]\times[0,2\pi]$, and apply periodic boundary
conditions in both directions. By discretizing the derivatives using
second-order central differences in space we can write this equation
in the discretized form
\begin{equation}
  \frac{\text{d} {\bf q}}{\text{d} t}={\bf A}{\bf q}+{\bf f}\sin(\omega t),
\end{equation}
where the vectors ${\bf q}$ and ${\bf f}$ denote the state and forcing
at each grid point, respectively. The state matrix ${\bf A}$
encapsulates the homogeneous component of the original system together
with the boundary conditions.

To formulate a direct-adjoint looping problem we consider choosing a
force ${\bf f}={\bf f}_{\textrm{true}}$ and evolving the system over
the time interval $[0,T]$ starting from a zero initial condition. This
results in the solution ${\bf q}_{\textrm{true}}$. We seek to recreate
the forcing ${\bf f}_{\textrm{true}}$ using only the observed
solution, ${\bf q}_{\textrm{true}}$. To this end, we let the cost
functional be
\begin{equation}
  \mathcal{J}=\frac{1}{T}\int_0^T \|{\bf q}-{\bf
    q}_\textrm{true}\|^2\,\textrm{d}t,
  \label{equ:cost}
\end{equation}
giving the adjoint equation
\begin{equation}
  \frac{\text{d} {\bf q}^\dagger(\tau)}{\text{d} \tau}={\bf A}^T{\bf
    q}^\dagger(\tau)+2\left({\bf q}(\tau)-{\bf
    q}_{\textrm{true}}(\tau)\right),
\end{equation}
where we have made the substitution $\tau = T-t$ to obtain an equation
we integrate forward in the time variable $\tau$. Once integrated, we
can use the adjoint solution to evaluate the gradient of the cost
functional with respect to the current forcing ${\bf f}$ as
\begin{equation}
  \frac{\partial \mathcal{L}}{\partial {\bf f}}=\frac{1}{T}\int_0^T
       {\bf q}^\dagger \sin(\omega t)\,\textrm{d}t.
       \label{equ:update}
\end{equation}
This gradient can then be used as part of a gradient-based
optimization routine to minimize the cost functional and hence find
${\bf f}_{\textrm{true}}.$

As mentioned previously, the algorithms developed perform better for
stiffer systems. One way of introducing stiffness into our
advection-diffusion problem is through the diffusion parameter
$D$. This can be observed by studying the eigenvalues of ${\bf
  A}$. Increasing $D$ stretches the spectrum along the negative real
axis and thus makes our system progressively more stiff. For our
implementation we set $a=1$, $\omega=1$, $T=10$, $ f_{\textrm{true}}=\sin(x)\sin(y)$, vary
$D\in\{0.01,0.1,1,10\}$ and determine the average time to perform one
direct-adjoint loop over $3$ runs for a range of processor numbers. An initial guess of f=1 is used. In
this manner, we can observe to what extent the theoretical scalings
agree with the scalings found from the numerical experiments. We choose the final time sufficiently large so that the size of the time-steps taken by the solvers is small compared to the size of the time partitions, whilst, at the same time, small enough so that the total runtime is kept reasonable. The choice of $f_{\textrm{true}}$ and the initial guess are found to not affect the scalings obtained.

The code is written in \textit{python}, and a fourth-order Runge-Kutta scheme
(with adaptive time-steps based on a fifth-order Runge-Kutta error
approximation) is used for the inhomogeneous integrations. This
time-stepper is implemented via the inbuilt \textit{RK45} method of
the \textit{solve\_ivp} function contained in the \textit{integrate}
submodule of the \textit{scipy} library. For homogeneous integrations
we use our implementation of the ReLPM as described
by~\citet{realLeja}. The relative error tolerance for each of these
integrators is kept the same at $\textit{rtol}=10^{-3}.$

\begin{figure}[h]
  \center
  \includegraphics[width=0.7\textwidth]{ 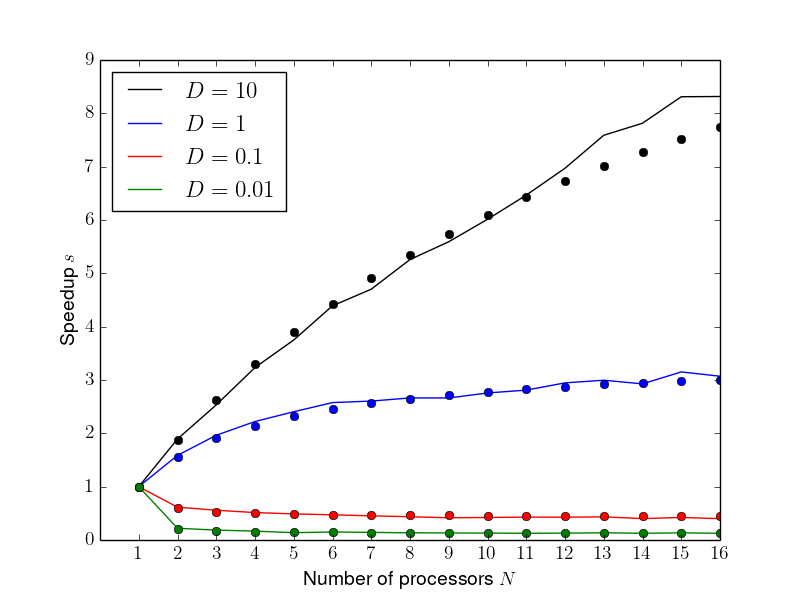}
  \caption{A figure showing the linear speedup for
    $D\in\{0.01,0.1,1,10\}$. The solid circles show the predicted
    values by the theoretical scaling formula (\ref{equ:LinSpeedup}).}
  \label{fig:speedupLin}
\end{figure}

Figure~\ref{fig:speedupLin} shows the speedup obtained for the linear
algorithm. As predicted, when $D$ increases the speedup becomes
larger, highlighting the observation that the homogeneous equation
needs to be solved significantly faster than the inhomogeneous
one. Furthermore, due to the extra cost in forming the matrix
exponential which in turn causes the homogeneous integration to be
slower than the inhomogeneous one, $D=0.01$ and $D=0.1$ report a worse
performance for the parallel algorithm when compared to the serial
approach. For $D=1$ and $D=10$ the stiffness introduced into the
system is circumvented through the exponential time-stepping and thus,
for the homogeneous equation, we are able to take far larger time
steps than the inhomogeneous integrator, off-setting the cost of
forming the matrix exponential. This observation can be confirmed by using the values for $\tau_H$ and $\tau_I$ from the serial running. We see that for $D=0.01$ and $D=0.1$ switching to the exponential integration is around 20 and 6.67 times slower, respectively. Whereas, for $D=1$ the exponential integration is 1.5 times faster, and for $D=10$ an even larger speedup of 4.3 is found.

The figure also shows good agreement between the numerically obtained
scaling and the one predicted by the theoretical value. This shows
that the effectiveness of a parallel-in-time approach for a specific
problem can be assessed simply by performing a short simulation to
determine the parameters $\tau_I$, $\tau^\dagger_I$, $\tau^\dagger_H$
and , $\tau^\dagger_H$. Furthermore, by assessing the theoretical
scaling the number of processors can best be chosen. In all cases, we
can clearly see the scaling approaching its maximum value as the
number of processors is increased. Initially, there are larger gains
in speedup, followed by progressively smaller gains as more processors
are included. For example, taking the $D=1$ case for instance: using
six processors we can expect a speedup of roughly $2.5,$ while
increasing the number of processors only provides a speed up of
approximately $3$ -- a relatively small speedup for more than double
the number of processors. Hence, the user should initially consider
the theoretical scaling and choose the number of processors that best
balances speedup gains versus expended resources.

\begin{table}
\centering
\begin{tabular}{rclclclcl}
 \hline 
 $N$ &\phantom{1}& $D=0.01$ &\phantom{123} & $D=0.1$ &\phantom{123} &  $D=1$ &\phantom{123} & $D=10$\\
\hline
 1 && 0.046 && 0.0034 && 0.0011 && 0.00061 \\ 
 2 && 0.058 && 0.0036 && 0.00087 && 0.0019 \\ 
 4 && 0.065 && 0.0088 && 0.0012 && 0.0040 \\ 
 8 && 0.074 && 0.020 && 0.0019 && 0.0017\\ 
 16 && 0.089 && 0.0072 && 0.0015 && 0.0040 \\ 
 \hline
\end{tabular}
\caption{Relative errors for the gradient obtained with a linear governing equation. The relative errors are calculated with respect to a more accurate gradient obtained in series.}
\label{table:linErrs}
\end{table}

To conclude our assessment of the linear algorithm we examine the error properties of performing a direct-adjoint loop in parallel. An accurate value for the gradient using equation (\ref{equ:update}) is first obtained in series using a higher error tolerance of $\textit{rtol}=10^{-8}$. Table \ref{table:linErrs} shows the relative error of the gradients computed with our lower error tolerance for different numbers of processors. We see that in general the error obtained by increasing the number of processors is not significantly different than those obtained in series. Slightly larger errors are observed when the number of processors is increased; this is in line with the original {\tt Paraexp} algorithm \citep{GanderParaExp}. The larger errors obtained for $D=0.01$ can be attributed to the solvers taking larger timesteps for this equation, leading to a reduced accuracy for both the interpolation of the direct equation and the computation of the gradient via numerical integration. This suggests that care must be taken to ensure that the time resolution is sufficient for computing these quantities.

\subsection{Non-linear governing equation}

To test the two algorithms for non-linear direct equations we consider
the harmonically forced viscous Burgers' equation

\begin{eqnarray}
  \frac{\partial U}{\partial t}&=&-U\frac{\partial U}{\partial
    x}-V\frac{\partial U}{\partial y} +D\nabla^2 U+ f_U\sin(\omega
  t),\\
  \frac{\partial V}{\partial t}&=&-U\frac{\partial V}{\partial
    x}-V\frac{\partial V}{\partial y} +D\nabla^2 V+ f_V\sin(\omega t),
\end{eqnarray}
for the state ${\bf q}=(U,V)^T$ on the same domain as the
advection-diffusion equation of the previous section. Once again, we
discretise using second-order centered finite differences and,
together with periodic boundary conditions, arrive at
\begin{equation}
  \frac{\text{d}{\bf q}}{\text{d} t}={\mathcal N}({\bf q})+{\bf f}\sin(\omega t).
\end{equation}
The equations linearised about the direct solution ${\bf q}=({\bf
  U},{\bf V })^T$ are
\begin{eqnarray}
  \frac{\partial u'}{\partial t}&=&-U\frac{\partial u'}{\partial
    x}-V\frac{\partial u'}{\partial y}-u'\frac{\partial U}{\partial
    x}-v'\frac{\partial U}{\partial y} +D\nabla^2 u',\\
  \frac{\partial v'}{\partial t}&=&-U\frac{\partial v'}{\partial
    x}-V\frac{\partial v'}{\partial y}-u'\frac{\partial V}{\partial
    x}-v'\frac{\partial V}{\partial y} +D\nabla^2 v',
\end{eqnarray}
for the linearised state ${\bf q}'=({\bf u}',{\bf v}')^T$. These
equations are also discretized using second-order centered finite
differences with periodic boundary conditions. This leads to the
linearized equations in the form
\begin{equation}
  \frac{\text{d}{\bf q}'}{\text{d}t}={\bf A}(U(t),V(t)){\bf q}', 
\end{equation}
with the corresponding adjoint equation 
\begin{equation}
  \frac{\text{d}{\bf q}^\dagger}{\text{d}\tau}={\bf
    A}^T(U(\tau),V(\tau)){\bf q}^\dagger+2({\bf q}(\tau)-{\bf
    q}_\textrm{true}(\tau)).
\end{equation}
We have made explicit the dependence of the state matrix ${\bf A}$ on
the direct solution and have, as before, introduced the time variable
$\tau=T-t$ to integrate the adjoint equations forwards in the time
variable $\tau$. The parameters chosen are $T=1$, $\omega=1$ and $ f_{U,\textrm{true}}=f_{V,\textrm{true}}=\sin(x)\sin(y)$ with an initial guess of $f_U=f_V=1$. Again, the choice of ${\bf f}_{\textrm{true}}$ and the initial guess are found to not affect the results.

\subsubsection{Non-linear governing equation solved in parallel}

\begin{figure}[h]
  \center
  \includegraphics[width=0.7\textwidth]{ 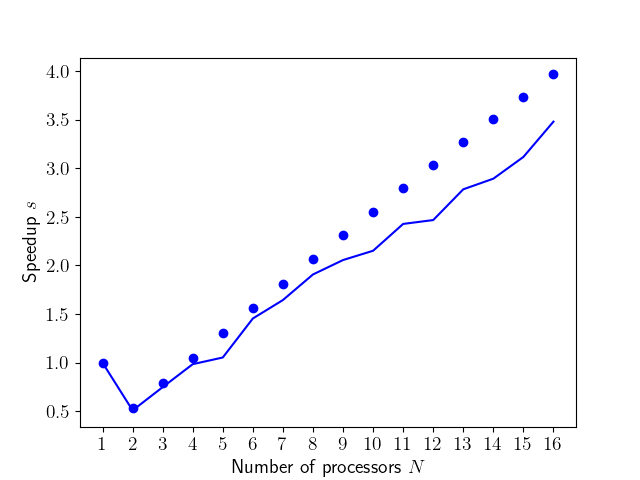}
  \caption{Non-linear speedup with the iterative procedure for
    $D=1$. The solid circles show the predicted values by the
    theoretical scaling formula~(\ref{equ:nonLinSpeedup}).}
  \label{fig:speedupNonLin}
\end{figure}

Now we turn to timing the direct-adjoint loop, i.e., we solve both the
non-linear governing equation and the adjoint in parallel as described
by algorithm~\ref{alg:nonlin}. Figure~\ref{fig:speedupNonLin} shows
the comparison between the numerically obtained speedup and the
speedup predicted by our theoretical
scaling~(\ref{equ:nonLinSpeedup}). As four iterations are used to
converge the direct solution, we set $K=4$ when computing the
theoretical scaling.

We note again that there is good agreement between the numerically
obtained speedup and that predicted by our theoretical
arguments. Initially the speedup is less than one due to the cost of
performing four iterations negating any possible speedup. However, for
$N>4$ this cost is overcome, and we obtain a marked speedup.

\subsubsection{Non-linear governing equation solved in series}
\label{sec:resultsHyb}

\begin{figure}[h]
  \center
  \includegraphics[width=0.7\textwidth]{ 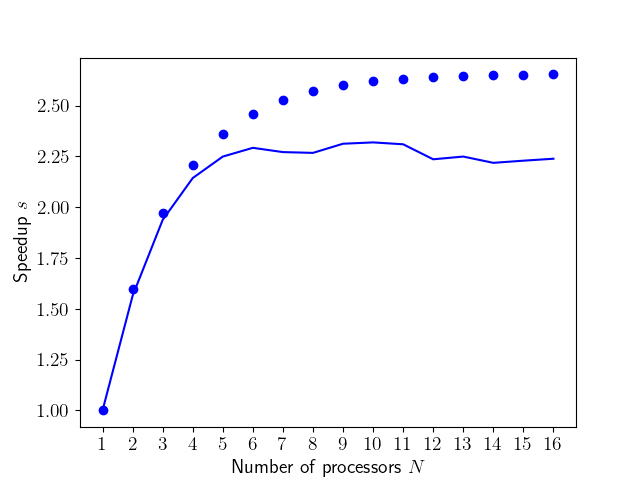}
  \caption{Non-linear speedup with the hybrid procedure for $D=1$. The
    solid circles show the predicted values by the theoretical scaling
    formula (\ref{equ:hybridSpeedup}).}
  \label{fig:speedupHybnonLin}
\end{figure}

Next, we time the same direct-adjoint loop, but using
algorithm~\ref{alg:hyb} in which the direct equation is solved in
series. By performing short direct and adjoint simulations we obtain
the value $k=2.11$ for calculating the size of the time partitions
using equation~(\ref{equ:nonEquiP}). The speedup obtained, compared to
the theoretical value, is shown in figure~\ref{fig:speedupHybnonLin}.

The figure shows that the speedup converges to its theoretical maximum
much faster than in the previous two algorithms. Indeed, about four
processors are needed for the speedup to obtain its theoretical
maximum. We also see a larger discrepancy in the obtained maximum
efficiency gain and the one predicted by our scaling equation. This
can mainly be attributed to the fact that the scaling argument assumes
that, as the width of the time partition gets smaller, the
inhomogeneous and homogeneous equations can be solved increasingly
faster. However, for our non-equispaced time partitioning, there comes
a point in which the time partition becomes smaller than the minimum
time-step needed to solve the equation. Hence, decreasing the width of
the time partition beyond this point does not yield any speedup in
solving the equation.

\begin{table}
\centering
\begin{tabular}{rclcl} 
 \hline 
 $N$ &\phantom{123} & iterative &\phantom{123}& hybrid \\
\hline
 1  && 0.0030 && 0.0030  \\ 
 2  && 0.0044 && 0.0027   \\ 
 4  && 0.0030 && 0.0033  \\ 
 8  && 0.0057 && 0.0028  \\ 
 16 && 0.0097 && 0.0029  \\ 
 \hline
\end{tabular}
\caption{Relative errors for the gradient obtained with a non-linear governing equation. The relative errors are calculated with respect to a more accurate gradient obtained in series with a relative tolerance of $\textit{rtol}=10^{-8}$.}
\label{table:nonLinErrs}
\end{table}

Similarly to the linear case, table \ref{table:nonLinErrs} shows the errors for the gradients obtained for the non-linear equation. Again, we see that increasing the number of processors does not change the error appreciably compared to the serial case. Interestingly, the error obtained with the hybrid algorithm stays more constant than for the case of solving the non-linear equation in parallel. This is to be expected as there is an increased source of error for non-linear equations stemming from solving the equation iteratively which is off-set by using the hybrid algorithm.

\subsection{Checkpointing}

\begin{figure}[h]
  \center
  \includegraphics[width=0.7\textwidth]{ 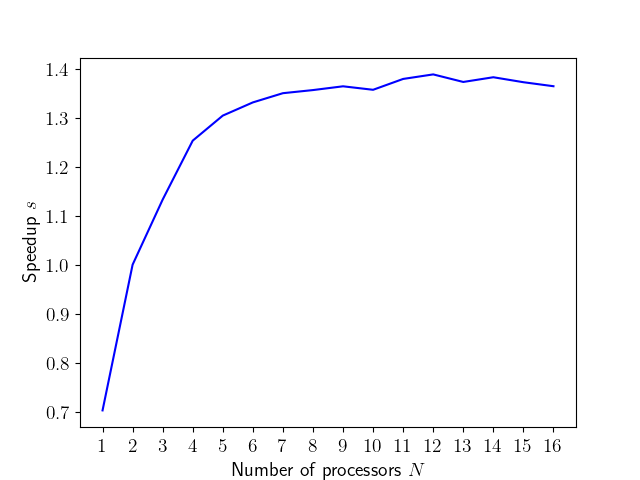}
  \caption{Non-linear speedup with the hybrid procedure for $D=1$ when
    five checkpoints are used.}
  \label{fig:speedupHybnonLinCheck}
\end{figure}

Lastly, we turn our attention to the inclusion of a checkpointing
scheme in performing a direct-adjoint loop. We consider the same
problem as in section~\ref{sec:resultsHyb}, but we include five
checkpoints using the method outlined in
section~\ref{sec:checkpointing}. The speedup obtained is shown in
figure~\ref{fig:speedupHybnonLinCheck}.

Similarly to the hybrid example with no checkpoints, we obtain the
maximum speedup at about four processors. The maximum speedup has
decreased, but this is to be expected for the hybrid algorithm as no
speedup is gained in solving the forward equation by itself. Only
direct-adjoint loops are accelerated. The fact that a speedup is still
obtained with the inclusion of checkpoints demonstrates the
compatibility of our developed algorithms with checkpointing
schemes. For the cases of a linear governing equation, and a
non-linear direct equation solved in parallel, we can expect less of a
penalty for using a checkpointing scheme as, in this case, the forward
equations are also accelerated.

\section{Conclusions}

In this article we have presented three separate algorithms that aim
to extend the {\tt{Paraexp}} algorithm of~\citet{GanderParaExp} to
adjoint-looping studies. Theoretical speed-ups have been derived for
all cases and verified through a set of numerical experiments. The
verifications of the theoretical predictions indicate that studies
utilizing these approaches can {\it{a-priori}} assess the achievable
speedups as well as the necessary processors required. All algorithms
sped up the calculation of the relevant system in line with the
predictions we presented, signifying that the adjoint {\tt{Paraexp}}
algorithm, and its non-linear extension, can be used to accelerate
direct-adjoint optimization studies. This is particularly pertinent
due to the linear, time-varying nature of the adjoint which lends
itself naturally to this parallelization approach, and thus can always
be parallelized in time.

We first considered the extension of the {\tt{Paraexp}} algorithm to
direct-adjoint loops that stem from a linear governing equation. As
the adjoint equation is also a linear equation this can be achieved
simply by integrating the adjoint equation using the {\tt{Paraexp}}
algorithm. By keeping solutions on each time partition local to each
processor we were able to minimize the processor communications
required for direct-adjoint looping. The theoretical scaling of the
algorithm was derived, showing that this algorithm has a theoretical
maximum that can be obtained as the number of processors is
increased. The linear algorithm was demonstrated on a two-dimensional
advection-diffusion equation, where good agreement with theoretical
scalings was obtained. These numerical experiments, as well as the
theoretical scalings, highlight an important consideration that must
be made when using the {\tt{Paraexp}} algorithm, namely that the
stiffness of the governing equations should be taken into
account. Only when the homogeneous component can be solved
significantly faster than the inhomogeneous one, can speedups be
obtained. Keeping this caveat in mind, we expect this approach to be
particularly applicable to systems governing by a stiff (multi-scale)
dynamics.

Further to the linear case, we also considered how a parallel-in-time
approach can be used for non-linear direct equations. This led to the
development of two algorithms; an iterative algorithm in which the
direct equation is solved in parallel, and a hybrid algorithm in which
the direct equation is solved in series. In both cases the adjoint
equation is integrated in parallel using the {\tt{Paraexp}}
algorithm. The scalings for the iterative approach showed that due to
iterations being carried out, care must be taken to ascertain whether
a speedup is possible for a given number of processors as the cost of
iterations can readily outweigh the cost of solving the equations in
parallel. Also, similarly to the linear case, there is a maximum
speedup possible. Therefore, an initial analysis of the governing
equations may be necessary to ascertain if this approach is
applicable and viable.

As in the iterative non-linear algorithm the iterations may not
converge, or too many processors may be needed in order to obtain a
speedup, we also proposed a hybrid approach in which the direct
equation is solved in series whilst the adjoint is solved in
parallel. By using a non-equidistant time-partitioning and overlapping
the direct and adjoint solves, a speedup can be obtained. The
time-partitioning is chosen so that all the homogeneous adjoint solves
begin at the same time. In this way, the theoretical scaling shows
that, as the number of processors increases, we obtain a maximum
speedup where the inhomogeneous adjoint integration is effectively
replaced by a homogeneous one. The scaling and subsequent numerical
verification demonstrated that in this approach the maximum speedup,
although lower than the other algorithms, is obtained with a smaller
number of processors. Again, we note that for the non-linear
algorithms the speedups rely on the homogeneous equations being able
to be solved more efficiently than their inhomogeneous counterparts,
for example as in the case of a stiff system.

Lastly, we demonstrated that all developed algorithms are easily
implementable with checkpointing regimes. Hence, we believe that a {\tt Paraexp} based parallel-in-time approach can provide a viable option for
accelerating direct-adjoint loops stemming from both linear and
non-linear governing equations, with the theoretical scalings
providing an easy `{\it{a-priori}} check' of the speedups
available. We emphasise that this approach is particularly aimed at
studies that have exhausted spatial parallelization but still have
additional computational power available. In this case, a significant
impact can be made by parallelizing in time, adding valuable
efficiency to these optimization studies. The optimal overall balance
of spatial versus temporal parallelism for a given number of
processors and a given computer architecture is an interesting
extension of our study and will be addressed in a future effort.

\section*{Acknowledgments}
The authors wish to acknowledge the EPSRC and Roth PhD scholarships on
which this research was conducted. We are also grateful to the HPC service at Imperial College London for providing the computational resources used for this study.

\bibliography{paraAdjResponseNoColours}

\begin{thebibliography}{50}
\providecommand{\natexlab}[1]{#1}
\providecommand{\url}[1]{\texttt{#1}}
\expandafter\ifx\csname urlstyle\endcsname\relax
  \providecommand{\doi}[1]{doi: #1}\else
  \providecommand{\doi}{doi: \begingroup \urlstyle{rm}\Url}\fi

\bibitem[Bal and Maday(2002)]{Guillaume2002}
G.~Bal and Y.~Maday.
\newblock A ``parareal'' time discretization for non-linear pde's with
  application to the pricing of an {A}merican put.
\newblock In \emph{Recent Developments in Domain Decomposition Methods}, pages
  189--202, Berlin, Heidelberg, 2002. Springer Berlin Heidelberg.

\bibitem[Bergamaschi et~al.(2006)Bergamaschi, Caliari, Mart{\'i}nez, and
  Vianello]{realLeja}
L.~Bergamaschi, M.~Caliari, A.~Mart{\'i}nez, and M.~Vianello.
\newblock Comparing {L}eja and {K}rylov approximations of large scale matrix
  exponentials.
\newblock In \emph{Computational Science -- ICCS 2006}, volume 3994 of
  \emph{Lect. Notes in Comp. Sci.}, pages 685--692, Berlin, Heidelberg, 2006.
  Springer-Verlag.

\bibitem[Caliari(2007)]{caliari2007}
M.~Caliari.
\newblock Accurate evaluation of divided differences for polynomial
  interpolation of exponential propagators.
\newblock \emph{Computing}, 80\penalty0 (2):\penalty0 189--201, 2007.

\bibitem[Clarke et~al.(2020{\natexlab{a}})Clarke, Davies, Ruprecht, and
  Tobias]{Clarke2020a}
A.~T. Clarke, C.~J. Davies, D.~Ruprecht, and S.~M. Tobias.
\newblock Parallel-in-time integration of kinematic dynamos.
\newblock \emph{Journal of Computational Physics: X}, 7:\penalty0 100057,
  2020{\natexlab{a}}.

\bibitem[Clarke et~al.(2020{\natexlab{b}})Clarke, Davies, Ruprecht, Tobias, and
  Oishi]{Clarke2020b}
A.~T. Clarke, C.~J. Davies, D.~Ruprecht, S.~M. Tobias, and J.~S. Oishi.
\newblock {Performance of parallel-in-time integration for Rayleigh
  B{\'{e}}nard convection}.
\newblock \emph{Computing and Visualization in Science}, 23\penalty0
  (1):\penalty0 10, 2020{\natexlab{b}}.

\bibitem[Cox and Matthews(2002)]{cox2002}
S.~M. Cox and P.~C. Matthews.
\newblock Exponential time differencing for stiff systems.
\newblock \emph{J. Comp. Phys.}, 176\penalty0 (2):\penalty0 430--455, 2002.

\bibitem[Eggl and Schmid(2018)]{Eggl2018}
M.~F. Eggl and P.~J. Schmid.
\newblock A gradient-based framework for maximizing mixing in binary fluids.
\newblock \emph{J. Comp. Phys.}, 368:\penalty0 131--153, 2018.

\bibitem[Eggl and Schmid(2020{\natexlab{a}})]{Eggl2019b}
M.~F. Eggl and P.~J. Schmid.
\newblock Mixing enhancement in binary fluids using optimised stirring
  strategies.
\newblock \emph{Journal of Fluid Mechanics}, 899:\penalty0 A24,
  2020{\natexlab{a}}.

\bibitem[Eggl and Schmid(2020{\natexlab{b}})]{eggl2019}
M.~F. Eggl and P.~J Schmid.
\newblock {Shape optimization of stirring rods for mixing binary fluids}.
\newblock \emph{IMA Journal of Applied Mathematics}, pages 762--789,
  2020{\natexlab{b}}.

\bibitem[Emmett and Minion(2012)]{EmmettMinion2012}
M.~Emmett and M.~L. Minion.
\newblock {Toward an Efficient Parallel in Time Method for Partial Differential
  Equations}.
\newblock \emph{Communications in Applied Mathematics and Computational
  Science}, 7:\penalty0 105--132, 2012.

\bibitem[Falgout et~al.(2014)Falgout, Friedhoff, Kolev, MacLachlan, and
  Schroder]{Falgout2014}
R.~D. Falgout, S.~Friedhoff, Tz.~V. Kolev, S.~P. MacLachlan, and J.~B.
  Schroder.
\newblock Parallel time integration with multigrid.
\newblock \emph{SIAM Journal on Scientific Computing}, 36\penalty0
  (6):\penalty0 C635--C661, 2014.

\bibitem[Farhat and Chandesris(2003)]{Farhat2003}
C.~Farhat and M.~Chandesris.
\newblock Time-decomposed parallel time-integrators: theory and feasibility
  studies for fluid, structure, and fluid–structure applications.
\newblock \emph{Int. J. Num. Meth. Eng.}, 58\penalty0 (9):\penalty0 1397--1434,
  2003.

\bibitem[Foures et~al.(2014)Foures, Caulfield, and Schmid]{Foures2014}
D.~P.~G. Foures, C.~P. Caulfield, and P.~J. Schmid.
\newblock Optimal mixing in two-dimensional plane {Poiseuille} flow at finite
  {P\'eclet} number.
\newblock \emph{J. Fluid Mech.}, 748:\penalty0 241--277, 2014.

\bibitem[Friedhoff et~al.(2013)Friedhoff, Falgout, Kolev, MacLachlan, and
  Schroder]{FriedhoffEtAl2013}
S.~Friedhoff, R.~D. Falgout, T.~V. Kolev, S.~P. MacLachlan, and J.~B. Schroder.
\newblock {A Multigrid-in-Time Algorithm for Solving Evolution Equations in
  Parallel}.
\newblock In \emph{{Presented at: Sixteenth Copper Mountain Conference on
  Multigrid Methods, Copper Mountain, CO, United States, Mar 17 - Mar 22,
  2013}}, 2013.

\bibitem[Gander(1996)]{Gander1996}
M.~J. Gander.
\newblock Overlapping schwarz for linear and nonlinear parabolic problems.
\newblock In \emph{9th International Conference on Domain Decomposition
  Methods}, pages 97--104, 1996.

\bibitem[Gander(2014)]{gander50years}
M.~J. Gander.
\newblock \emph{50 Years of Time Parallel Time Integration}, volume~9 of
  \emph{Contributions in Mathematical and Computational Sciences}, pages
  69--113.
\newblock Springer International Publishing, Cham, 2014.

\bibitem[Gander and G\"{u}ttel(2013)]{GanderParaExp}
M.~J. Gander and S.~G\"{u}ttel.
\newblock {PARAEXP}: A parallel integrator for linear initial-value problems.
\newblock \emph{SIAM J. Sci. Comp.}, 35\penalty0 (2):\penalty0 C123--C142,
  2013.

\bibitem[Gander and Neum{\"{u}}ller(2016)]{Gander2016}
M.~J. Gander and M.~Neum{\"{u}}ller.
\newblock Analysis of a new space-time parallel multigrid algorithm for
  parabolic problems.
\newblock \emph{SIAM Journal on Scientific Computing}, 38\penalty0
  (4):\penalty0 A2173--A2208, 2016.

\bibitem[Gander et~al.(2018)Gander, G{\"u}ttel, and Petcu]{Gander2017ANP}
M.~J. Gander, S.~G{\"u}ttel, and M.~Petcu.
\newblock A nonlinear paraexp algorithm.
\newblock In \emph{Domain Decomposition Methods in Science and Engineering
  XXIV}, pages 261--270, Cham, 2018. Springer International Publishing.

\bibitem[Gander et~al.(2020{\natexlab{a}})Gander, Kwok, and
  Salomon]{Gander2020}
M.~J. Gander, F.~Kwok, and J.~Salomon.
\newblock Paraopt: A parareal algorithm for optimality systems.
\newblock \emph{SIAM Journal on Scientific Computing}, 42\penalty0
  (5):\penalty0 A2773--A2802, 2020{\natexlab{a}}.

\bibitem[Gander et~al.(2020{\natexlab{b}})Gander, Liu, Wu, Yue, and
  Zhou]{GanderEtAl2020}
M.~J. Gander, J.~Liu, S.-L. Wu, X.~Yue, and T.~Zhou.
\newblock Paradiag: Parallel-in-time algorithms based on the diagonalization
  technique.
\newblock 2020{\natexlab{b}}.
\newblock URL \url{http://arxiv.org/abs/2005.09158}.

\bibitem[Giladi and Keller(2002)]{Giladi2002}
E.~Giladi and H.~B. Keller.
\newblock {Space-time domain decomposition for parabolic problems}.
\newblock \emph{Numerische Mathematik}, 93\penalty0 (2):\penalty0 279--313,
  2002.

\bibitem[G{\"o}tschel and Minion(2018)]{Gotschel2018}
S.~G{\"o}tschel and M.~L. Minion.
\newblock Parallel-in-time for parabolic optimal control problems using pfasst.
\newblock In \emph{Domain Decomposition Methods in Science and Engineering
  XXIV}, pages 363--371, Cham, 2018. Springer International Publishing.

\bibitem[G{\"o}tschel and Minion(2019)]{Gotschel2019}
S.~G{\"o}tschel and M.~L. Minion.
\newblock An efficient parallel-in-time method for optimization with parabolic
  pdes.
\newblock \emph{SIAM Journal on Scientific Computing}, 41\penalty0
  (6):\penalty0 C603--C626, 2019.

\bibitem[Griewank and Walther(2000)]{revolve}
A.~Griewank and A.~Walther.
\newblock Algorithm 799: Revolve: An implementation of checkpointing for the
  reverse or adjoint mode of computational differentiation.
\newblock \emph{ACM Trans. Math. Softw.}, 26:\penalty0 19--45, 2000.

\bibitem[G{\"{u}}nther et~al.(2018)G{\"{u}}nther, Gauger, and
  Schroder]{Gunther2018}
S.~G{\"{u}}nther, N.~R. Gauger, and J.~B. Schroder.
\newblock {A non-intrusive parallel-in-time adjoint solver with the XBraid
  library}.
\newblock \emph{Computing and Visualization in Science}, 19\penalty0
  (3):\penalty0 85--95, 2018.

\bibitem[G{\"u}nther et~al.(2019)G{\"u}nther, Gauger, and
  Schroder]{Gunther2019}
S.~G{\"u}nther, N.~R. Gauger, and J.~B. Schroder.
\newblock A non-intrusive parallel-in-time approach for simultaneous
  optimization with unsteady pdes.
\newblock \emph{Optimization Methods and Software}, 34\penalty0 (6):\penalty0
  1306--1321, 2019.

\bibitem[G{\"u}ttel et~al.(2020)G{\"u}ttel, Kressner, and
  Lund]{matrixFunctions}
S.~G{\"u}ttel, D.~Kressner, and K.~Lund.
\newblock Limited-memory polynomial methods for large-scale matrix functions.
\newblock \emph{GAMM-Mitteilungen}, 43\penalty0 (3):\penalty0 e202000019, 2020.

\bibitem[Horton and Vandewalle(1995)]{HortonVandewalle1995}
G.~Horton and S.~Vandewalle.
\newblock {A Space-Time Multigrid Method for Parabolic Partial Differential
  Equations}.
\newblock \emph{SIAM Journal on Scientific Computing}, 16\penalty0
  (4):\penalty0 848--864, 1995.

\bibitem[Jameson(1988)]{Jameson1988}
A.~Jameson.
\newblock Aerodynamic design via control theory.
\newblock \emph{J. Sci. Comp.}, 3\penalty0 (3):\penalty0 233--260, 1988.

\bibitem[Kallala et~al.(2019)Kallala, Vay, and Vincenti]{kallala2019}
H.~Kallala, J.-L. Vay, and H.~Vincenti.
\newblock A generalized massively parallel ultra-high order {FFT}-based
  {Maxwell} solver.
\newblock \emph{Comp. Phys. Comm.}, 244:\penalty0 25--34, 2019.

\bibitem[Kassam and Trefethen(2005)]{kassam2005}
A.-K Kassam and L.~N. Trefethen.
\newblock Fourth-order time-stepping for stiff {PDEs}.
\newblock \emph{SIAM J. Sci. Comp.}, 26\penalty0 (4):\penalty0 1214--1233,
  2005.

\bibitem[Kooij et~al.(2017)Kooij, Botchev, and Geurts]{ParaExpNL}
G.~L. Kooij, M.~A. Botchev, and B.~J. Geurts.
\newblock A block {Krylov} subspace implementation of the time-parallel
  {Paraexp} method and its extension for nonlinear partial differential
  equations.
\newblock \emph{J. Comp. Appl. Math.}, 316:\penalty0 229--246, 2017.

\bibitem[Kwok(2014)]{Kwok2014}
F.~Kwok.
\newblock Neumann--neumann waveform relaxation for the time-dependent heat
  equation.
\newblock In \emph{Domain Decomposition Methods in Science and Engineering
  XXI}, pages 189--198, Cham, 2014. Springer International Publishing.

\bibitem[Laizet and Vassilicos(2010)]{Laizet2010}
S.~Laizet and J.~C. Vassilicos.
\newblock Direct numerical simulation of fractal-generated turbulence.
\newblock In \emph{Direct and Large-Eddy Simulation VII}, pages 17--23,
  Dordrecht, 2010. Springer Netherlands.

\bibitem[{Lelarasmee} et~al.(1982){Lelarasmee}, {Ruehli}, and
  {Sangiovanni-Vincentelli}]{Lelarasmee1982}
E.~{Lelarasmee}, A.~E. {Ruehli}, and A.~L. {Sangiovanni-Vincentelli}.
\newblock The waveform relaxation method for time-domain analysis of large
  scale integrated circuits.
\newblock \emph{IEEE Transactions on Computer-Aided Design of Integrated
  Circuits and Systems}, 1\penalty0 (3):\penalty0 131--145, 1982.

\bibitem[Lions et~al.(2001)Lions, Maday, and Turinici]{lions2001}
J.-L. Lions, Y.~Maday, and G.~Turinici.
\newblock R{\'e}solution d'edp par un sch{\'e}ma en temps parar{\'e}el.
\newblock \emph{Compt. R. de l'Acad. Sci. I -- Math.}, 332\penalty0
  (7):\penalty0 661--668, 2001.

\bibitem[Maday and R{\o}nquist(2008)]{MADAY2008113}
Y.~Maday and E.~M. R{\o}nquist.
\newblock Parallelization in time through tensor-product space-time solvers.
\newblock \emph{Comptes Rendus Mathematique}, 346\penalty0 (1):\penalty0
  113--118, 2008.

\bibitem[Maday and Turinici(2002)]{Maday2002}
Y.~Maday and G.~Turinici.
\newblock A parareal in time procedure for the control of partial differential
  equations.
\newblock \emph{Comptes Rendus Mathematique}, 335\penalty0 (4):\penalty0
  387--392, 2002.

\bibitem[Maday and Turinici(2003)]{Maday2003}
Y.~Maday and G.~Turinici.
\newblock Parallel in time algorithms for quantum control: Parareal time
  discretization scheme.
\newblock \emph{International Journal of Quantum Chemistry}, 93\penalty0
  (3):\penalty0 223--228, 2003.

\bibitem[Mandal(2014)]{Mandal2014}
B.~C. Mandal.
\newblock A time-dependent dirichlet-neumann method for the heat equation.
\newblock In \emph{Domain Decomposition Methods in Science and Engineering
  XXI}, pages 467--475, Cham, 2014. Springer International Publishing.

\bibitem[Marcotte and Caulfield(2018)]{Marcotte2018}
F.~Marcotte and C.~P. Caulfield.
\newblock Optimal mixing in two-dimensional stratified plane {Poiseuille} flow
  at finite {P\'eclet} and {Richardson} numbers.
\newblock \emph{J. Fluid Mech.}, 853:\penalty0 359--385, 2018.

\bibitem[Nievergelt(1964)]{Nievergelt1964}
J.~Nievergelt.
\newblock Parallel methods for integrating ordinary differential equations.
\newblock \emph{Commun. ACM}, 7\penalty0 (12):\penalty0 731--733, 1964.

\bibitem[Ong and Schroder(2020)]{Ong2020}
B.~W. Ong and J.~B. Schroder.
\newblock {Applications of time parallelization}.
\newblock \emph{Computing and Visualization in Science}, 23\penalty0
  (1):\penalty0 11, 2020.

\bibitem[Pekurovsky(2012)]{pekurovsky2012}
D.~Pekurovsky.
\newblock {P3DFFT}: A framework for parallel computations of {Fourier}
  transforms in three dimensions.
\newblock \emph{SIAM J. Sci. Comp.}, 34\penalty0 (4):\penalty0 C192--C209,
  2012.

\bibitem[Pringle and Kerswell(2010)]{KerswellTurb}
C.~T. Pringle and R.~Kerswell.
\newblock Using nonlinear transient growth to construct the minimal seed for
  shear flow turbulence.
\newblock \emph{Phys. Rev. Lett.}, 105:\penalty0 154502, 2010.

\bibitem[Qadri et~al.(2016)Qadri, Magri, Ihme, and Schmid]{ubaidFlame}
U.~A. Qadri, L.~Magri, M.~Ihme, and P.~J. Schmid.
\newblock Optimal ignition placement in diffusion flames by nonlinear adjoint
  looping.
\newblock \emph{Ctr. Turb. Res., Proc. of the Summer Program}, 2016.

\bibitem[Saad(1992)]{krylovApprox}
Y.~Saad.
\newblock Analysis of some {K}rylov subspace approximations to the matrix
  exponential operator.
\newblock \emph{SIAM J. Num. Anal.}, 29\penalty0 (1):\penalty0 209--228, 1992.

\bibitem[Schulze et~al.(2009)Schulze, Schmid, and Sesterhenn]{schulze2009}
J.~C. Schulze, P.~J. Schmid, and J.~L. Sesterhenn.
\newblock Exponential time integration using {Krylov} subspaces.
\newblock \emph{Int. J. Num. Meth. Fluids}, 60\penalty0 (6):\penalty0 591--609,
  2009.

\bibitem[Skene(2019)]{skene}
C.~S. Skene.
\newblock \emph{Adjoint based analysis for swirling and reacting flows}.
\newblock PhD thesis, Imperial College London, 2019.

\end{thebibliography}

\end{document}